\documentclass[11pt]{article}

\usepackage{etex}
\reserveinserts{28}
\usepackage[all]{xy}

\usepackage{amsfonts}
\usepackage{amsthm}
\usepackage{graphicx}
\usepackage{enumerate}

\usepackage{mathrsfs}
\usepackage{bm}
\usepackage{cite}
\usepackage{amssymb,amsmath} 
\usepackage{makecell}
\usepackage{tikz}
\usepackage{pgf}
\usepackage{tikz}
\usetikzlibrary{patterns}
\usepackage{pgffor}
\usepackage{pgfcalendar}
\usepackage{pgfpages}
\usepackage{shuffle,yfonts}
\usepackage{mathtools}
\DeclareFontFamily{U}{shuffle}{}
\DeclareFontShape{U}{shuffle}{m}{n}{ <-8>shuffle7 <8->shuffle10}{}

\newcommand{\vU}{\mathsf {U}}

\newcommand{\baru}{{\bar u}}

\newcommand{\ga}{\alpha}
\newcommand{\gb}{\beta}
\newcommand{\gc}{\gamma}
\newcommand{\gd}{\delta}
\newcommand{\gG}{\Gamma}
\newcommand\bom{{\bar\omega}}
\newcommand\om{{\omega}}

\newcommand\eps{{\varepsilon}}

\newcommand\calC{{\mathcal{C}}}

\usetikzlibrary{arrows,shapes,chains}

\let\oldsection\section
\renewcommand\section{\setcounter{equation}{0}\oldsection}

\allowdisplaybreaks

\DeclareMathOperator{\ord}{ord}

\renewcommand{\gg}{{\epsilon}}

\def\R{\mathbb{R}}
\def\N{\mathbb{N}}\def\Z{\mathbb{Z}}
\def\Q{\mathbb{Q}}
\def\CC{\mathbb{C}}

\theoremstyle{plain}
\newtheorem{thm}{Theorem}

\newtheorem{cor}[thm]{Corollary}

\newtheorem{prop}[thm]{Proposition}
\theoremstyle{definition}

\begin{document}
 \title{\bf On $N$-Unital Functions}
\author{Jianqiang Zhao\thanks{Email: zhaoj@ihes.fr} \\
\ \\
\small Department of Mathematics, The Bishop's School, La Jolla, CA 92037, USA}

\date{}
\maketitle

\noindent{\bf Abstract.}
Let $N$ be a positive integer. We say a non-constant rational function $U(x)\in{\mathbb C}(x)$ is $N$-\emph{unital} if all the zeros and poles of both $U(x)$ and $1-U(x)$ are either 0 or $N$-th roots of unity. These functions are called admissible functions by Au in a recent paper arXiv:2007.03957 and used to study some central binomial series of Ap\'ery type via their iterated integral expression related to multiple polylogrithms and colored multiple zeta values. In this paper we determine the complete set of these functions for $N\le 4$ by elementary method, and briefly study some cases at level 5.




\bigskip

\section{Introduction}
The study of central binomial series of Ap\'ery type has gradually attracted the attention of a lot number theorist ever after
Ap\'ery \cite{Apery1978} successfully proved the irrationality of $\zeta(3)$ by using one of such series. Z.-W. Sun
posted many conjectured identities in \cite{Sun2020}. Oftentimes, these series are related to powers of $\pi$ of special values
of the Riemann zeta function. Recently, Au \cite{Au2020} showed that some of these series can expressed
using the colored multiple zeta values, i.e., special values of multiple polylogarithms at roots of unities.

Let $N$ be a positive integer.
We say a non-constant rational function $U(x)\in\CC(x)$ is $N$-unital if all the zeros and poles of both $U(x)$
and $1-U(x)$ are either 0 or $N$-th roots of unity. Let $\vU_N$ be the set of all $N$-unital functions.
For convenience, for any rational function $f$ we define its $S_3$ orbit by
$$ \langle f\rangle_6=\left\{ f,  1-f,\frac{1}{f},\frac{f}{f-1}, \frac{f-1}f, \frac{-1}{f-1}\right\}.
$$
Then it is not hard to show that $\vU_1=\langle x\rangle_6$. In \cite{Au2020}, Au used the name ``admissible''
for such functions. We changed the name to avoid confusion with admissibility of multiple zeta values
whose study is intimately related to these functions.

Au conjectured in \cite{Au2020} that  $\vU_N$ is finite for any $N\in\N$. As kindly pointed out to me by F.\ Brunault,
this in fact follows from a more general theorem on the $S$-unit equations of Mason \cite{Mason1983}, see \cite[Cor. 2.2]{Brunault2021}. He also graciously shared his MAGMA program he used to compute these functions which detected some missing 3- and 4-unital functions in a previous draft.

Au \cite{Au2020} made some conjectures on the set $\calC^N$ of possible values of such functions at
0 when $N=2$ and 4, which are crucial in his study of some  central binomial series of Ap\'ery type via their iterated integral expression related to multiple polylogrithms and colored multiple zeta values. In this paper we will determine $\vU_N$ precisely for $N\le 4$ so that we can compute  $\calC^2$ and  $\calC^4$ precisely. We also make a conjecture on $\calC^N$ for all $N$.

\section{Statements of the main results}
Au explained in \cite[Appendix B]{Au2020} how one might discover some $4$-unital functions,
and his list covers most of these function although it is not complete. In fact, his list contains
168 such functions while the true number should be 252.

For any positive integer $N$ we denote the set of $N$-th roots of unity by $\gG_N=\{\exp(2\pi i/j): j=1, \dots, N\}$.
First, Au observed that the orbit of $x$ under the octahedral symmetry $S_4$ of 24 elements and the $S_3$ symmetry of 6 elements
$\langle f\rangle_6$ produces 72 rational functions:
\begin{align*}
O(x):=& \bigcup_{\eps \in\gG_4}
\langle \eps x\rangle_6\cup
\left\langle \frac{ i(\eps x+1)}{\eps x-1}\right\rangle_{\hskip-4pt 6}\cup
\left\langle \frac{\eps x+1}{\eps x+i}\right\rangle_{\hskip-4pt 6}.
\end{align*}
Similarly, the orbit of $x^2$ has 36 rational functions:
\begin{align*}
O(x^2):=\bigcup_{\eps =\pm 1}
\langle \eps x^2\rangle_6\cup
\left\langle \eps \left(\frac{x-1}{x+1}\right)^2\right\rangle_{\hskip-4pt 6}\cup
\left\langle \eps \left(\frac{x-i}{x+i}\right)^2\right\rangle_{\hskip-4pt 6}.
\end{align*}
The  orbit of $x^4$ has 18 rational functions:
\begin{align*}
O(x^4):=
\langle x^4\rangle_6\cup
\left\langle \left(\frac{x+1}{x-1}\right)^4\right\rangle_{\hskip-4pt 6}\cup
\left\langle \left(\frac{x-i}{x+i}\right)^4\right\rangle_{\hskip-4pt 6}.
\end{align*}
Au also noticed the following two orbits: one with 36 elements
\begin{align*}
O\left(\frac{2x}{x^2+1}\right):=\bigcup_{\eps =\pm 1}
\left\langle \frac{2\eps x}{x^2+1}\right\rangle_{\hskip-4pt 6}\cup
\left\langle \frac{2\eps ix}{x^2-1}\right\rangle_{\hskip-4pt 6}\cup
\left\langle \frac{\eps(x^2-1)}{x^2+1}\right\rangle_{\hskip-4pt 6},
\end{align*}
and the other with 6 elements:  $\left\langle \frac{4x^2}{(x^2+1)^2}\right\rangle_{\hskip-2pt 6}$.

Now, it is not difficult to check that $\vU_4$ contains three more orbits
of 48, 24  and 12 elements (the non-trivial symmetries are from $C_4\times S_3$ where
$C_4$ is the cyclic symmetry determined by $i\to i^n$, $0\le n\le 4$):
\begin{align*}
O\left(\frac{2x}{x+1}\right):= &
\bigcup_{\eps =\pm 1}
\left\langle \frac{\eps(x-1)}{x+1}\right\rangle_{\hskip-4pt 6}\cup
\left\langle \frac{\eps(x-i)}{x+i}\right\rangle_{\hskip-4pt 6} \cup
\left\langle \frac{(1-i)x}{x+\eps}\right\rangle_{\hskip-4pt 6}\cup
\left\langle \frac{(1+i)x}{x+\eps}\right\rangle_{\hskip-4pt 6} ,\\
O\left(\frac{x(x-1)}{x^2+1}\right):= &
\left\langle \frac{x(x-1)}{x^2+1}\right\rangle_{\hskip-4pt 6}\cup
\left\langle \frac{x(x+1)}{x^2+1}\right\rangle_{\hskip-4pt 6}\cup
\left\langle \frac{x(x-i)}{x^2-1}\right\rangle_{\hskip-4pt 6}\cup
\left\langle \frac{x(x+i)}{x^2-1}\right\rangle_{\hskip-4pt 6} ,\\
O\left(\frac{2(1+i)x}{(x+1)(x+i)}\right):= &
\left\langle \frac{2(1+i)x}{(x+1)(x+i)}\right\rangle_{\hskip-4pt 6}\cup
\left\langle \frac{2(1-i)x}{(x+1)(x-i)}\right\rangle_{\hskip-4pt 6}.
\end{align*}
Graphically, we have the following explanation of the two orbits above:
\begin{center}
\includegraphics[height=1.5in]{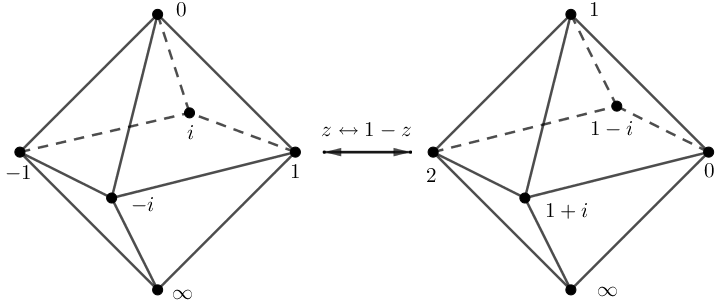}
\end{center}

\begin{thm}\label{thm:U4}
There are exactly $252$ $4$-unital functions which are given by
\begin{align*}
\vU_4= &\, O(x) \cup O(x^2) \cup O(x^4) \cup O\left(\frac{2x}{x^2+1}\right)\cup \left\langle \frac{4x^2}{(x^2+1)^2}\right\rangle_{\hskip-4pt 6} \\
&\, \cup O\left(\frac{2x}{x+1}\right) \cup O\left(\frac{x(x-1)}{x^2+1}\right)
\cup O\left(\frac{2(1+i)x}{(x+1)(x+i)}\right) .
\end{align*}
\end{thm}

\begin{cor}\label{cor:U2}
The set of all $2$-unital functions is given by
\begin{align*}
\vU_2&= \langle x\rangle\cup \langle -x\rangle\cup \langle x^2\rangle\cup \left\langle \frac{1+x}{2} \right\rangle_{\hskip-4pt 6}
\cup \left\langle \frac{1+x}{2x}\right\rangle_{\hskip-4pt 6}  \cup \left\langle \frac{(1+x)^2}{4x}\right\rangle_{\hskip-4pt 6} .
\end{align*}
\end{cor}
We list all the $2$- and $4$-unital functions explicitly in Appendix A.

Using the same approach we also find the complete set of $3$-unital functions and list them in Appendix B.
\begin{thm}\label{thm:U3}
Let $u=\exp(2\pi i/3)$. There are exactly $84$ $3$-unital functions which are given by
\begin{equation*}
\vU_3=\langle x^3\rangle_6 \cup
 \left\langle\frac{(1-\baru)x}{x-u}\right\rangle_{\hskip-4pt 6} \cup \bigcup_{\eps\in\gG_3} \langle \eps x\rangle_6 \cup \left\langle \frac{1-u}{\eps x-u}\right\rangle_{\hskip-4pt 6}  \cup\left\langle \frac{1-u\eps x}{1-\eps x}\right\rangle_{\hskip-4pt 6} \cup  \left\langle \frac{-3\eps x}{(\eps x-1)^2}\right\rangle_{\hskip-4pt 6}.
\end{equation*}
\end{thm}

For any positive integer $N$, we set
\begin{equation*}
\calC^N=\{U(0): U(x)\in \vU_N\}.
\end{equation*}

\begin{cor}
Let $u=\exp(2\pi i/3)$. We have
\begin{align*}
\calC^2=&\left\{0, \pm 1,\frac{1}{2}, 2,\infty \right\},  \quad
\calC^4=\left\{0, \pm 1,\frac{1}{2}, 2, \pm i, 1\pm i,\frac{1\pm i}{2},\infty  \right\}, \\
\calC^3=&\left\{0,  1, \pm u, \pm \baru,1-u, 1-\baru, \frac{1}{1-u}, \frac{1}{1-\baru}, \infty \right\}.
\end{align*}
\end{cor}

 From the explicit computation of this paper, we find the following result.
\begin{prop}
Let $N\le 4$ and $u_j=\exp(2j\pi i/N)$ for all $j=1,\dots, N-1$. Then
\begin{equation} \label{equ:ConjR0}
     \calC^N=\{0,1,\infty\} \cup
\left\{
  \begin{array}{ll}
  \displaystyle  \bigcup_{j=1}^{\frac{N}{2}} \langle \mu_j \rangle_6, & \hbox{if $N$ is even;} \\
  \displaystyle    \bigcup_{\eps=\pm1} \bigcup_{j=1}^{\frac{N-1}{2}} \langle \eps \mu_j \rangle_6, \qquad \  & \hbox{if $N$ is odd.}
  \end{array}
\right.
\end{equation}
Then $U(0)\in \calC^N$ for all $U(x)\in \vU_N$.
\end{prop}

Note that both $\langle 1 \rangle_6$ and  $\langle -1 \rangle_6$ each has only 3 elements. Moreover, the sets in the proposition may not be disjoint as shown when $N=3$.

Unfortunately, even for $N=5$ the set defined by the right-hand side of \eqref{equ:ConjR0} is too small to contain all values of 5-unital functions at 0. When $N=5$, setting $u=u_5=\exp(2\pi i/5)$  we find the following exceptional 5-unital functions:

\begin{align*}
f(x)=&\frac{(u^3 + u^2 + 2)x^2 + (u^4 + u^2)x - u^2 + u - 1}{x^2-(u^4 + u^2)x + u},\\
g(x)=&\frac{(u^2 + u + 1)x - u^2 - u - 1}{x-u^2}.   
\end{align*}
For example, one can easily check that 
$$
f(x)=\frac{(u^3 + u^2 + 2)(x- 1)(x-u)}{(x-u^2)(x-u^4)}, \quad
g(x)=\frac{(u^2 + u + 1)(x-1)}{x-u^2} 
\in\vU_5
$$
since
$$
1-f(x)=-\frac{(u^3 + u^2 + 1)(x-u^3)^2}{(x-u^2)(x-u^4)}, \quad
1-g(x)=\frac{-(u^2 + u)(x-u^4)}{x-u^2} 
\in\vU_5.
$$
But none of the following values are in the set on the right-hand side of \eqref{equ:ConjR0}:
\begin{alignat*}{6}
\langle f(0)\rangle_6=&\, \Big\{\frac{3-\sqrt{5}}2,  \frac{\sqrt{5}-1}2,  \frac{3+\sqrt{5}}2,
-\frac{1+\sqrt{5}}2, \frac{1-\sqrt{5}}2, \frac{\sqrt{5}+1}2\Big\}, \\
\langle g(0)\rangle_6=&\, \Big\{u+u^2, 1-u-u^2,-u-u^2, 1+u+u^2,-u^3-u^4 , 1+u^3+u^4  \Big\}.
\end{alignat*}

\section{Some general strategy}\label{sec:idea}
For any positive integer $N$ we denote the set of $N$-th roots of unity by $\gG_N=\{\exp(2\pi i/j): j=1, \dots, N\}$
and $\gG_N^0=\{0\}\cup \gG_N$. If $f\in\vU_N$, then we may assume $1-g=f$ for some $g\in\vU_N$. By considering
possible poles and zeros for $f$ and $g$, we see that $1-g=f$ can be transformed to the following form
\begin{equation}\label{equ:GenForm}
 \prod_{l\in L} (x-\mu_l)^{n_l}-C \prod_{k\in K} (x-\mu_k)^{\ga_k}=D\prod_{j\in J} (x-\mu_j)^{a_j}
\end{equation}
with $J \coprod K \coprod L\subset \gG_N^0$ where $\mu_0=0$ and all the exponents $n_j,\ga_j,a_j\in\N$.
Namely, we may assume $f$ and $g$ have exactly the same poles (counting multiplicities) but disjoint zero sets.

The following tricks (techniques?) will be used repeatedly in our proofs of the main theorems.
\begin{enumerate}
  \item[\bf{T1.}]  If $0\in K$ then we see that $D\in \pm \gG_N\subseteq \gG_{2N}$ by setting $x=0$. Similarly,
$0\in J \Rightarrow C\in \gG_{2N}$ while $0\in L \Rightarrow C/D\in \gG_{2N}$.

  \item[\bf{T2.}] Suppose $N=p$ is a prime. If $0,k\in K$ and $k\ne 0$ then
by setting $x=\mu_k$ we see that $\sum_{l\in L}  n_l=\sum_{j\in J} a_j$ by $p$-adic evaluation of
both sides of \eqref{equ:GenForm}. Similarly, $0,j\in J, j\ne 0\Rightarrow \sum_{l\in L}  n_l=\sum_{k\in K} \ga_k$
and $0,l\in L, l\ne 0\Rightarrow \sum_{k\in K} \ga_k=\sum_{j\in J} a_j$. We call this idea the \emph{$p$-adic argument.}

  \item[\bf{T3.}] Assume $N=p$ is a prime. Assume $K=\{0\}$ and $l\in L\ne\emptyset,j \in J\ne\emptyset$ then
by $x=0\Rightarrow D\in \gG_{2N}$. If $\sum_{j\in J} a_j<\ga_0$ then the highest degree must be $\ga_0$
on the LHS and thus $-C=D$ by comparing the leading coefficient. But then
$x=\mu_l\Rightarrow -\mu^{\ga_0}=\prod_{j\in J} (\mu_l-\mu_j)^{a_j}$ which is absurd by the $p$-adic evaluation.
Similarly, $\sum_{j\in L} k_j\ge \ga_0$ and therefore $\sum_{j\in L} k_j=\sum_{j\in J} a_j$.
We call this the \emph{single power $N$-adic argument.}

\item[\bf{T4.}] If $N=4$, we will have to modify the above argument by considering the $2$-divisibility of the terms after
evaluating $x$ at $\gG_4^0$. For any
$0\ne a\in \Q[i]$, define $\ord_2(a)=v$ if $|a|=2^v u$ where $u$ is a 2-unit, i.e., $u=p/q$ where $p,q$ are both odd integers.
We call this idea of considering $2$-divisibility  $2$-adic argument.
\end{enumerate}

We will mainly apply T2 and T3 when $N=3$. Note that the $p$-adic argument can imply impossibility
if one of $J,K,L$ is empty while the single power $p$-adic argument cannot since $J,K,L$ are always nonempty
when this argument applies.

In the proof of the main theorems, we often need to compares coefficients of polynomials.
For convenience, we use the shorthand $n$HD for ``the $n$th highest degree term''
for all $n\ge 2$. For example, 2HD stands for ``the second highest degree term''.

\section{Proof of Theorem \ref{thm:U4}: complete set of 4-unital functions}
Suppose $f= D x^a (x-1)^b (x+i)^c (x-i)^d (x+i)^e, g=C x^\ga (x-1)^\gb (x+i)^\gc (x-i)^\gd (x+i)^\gg \in\vU_4$ such that $1-g=f$.
We now break into different cases according the signs of the exponents. Observe that we can simplify our proof by using the fact
that $U(x)\in \vU_4 \Leftrightarrow U(-x)\in \vU_4 \Leftrightarrow U(ix)\in \vU_4\Leftrightarrow U(-ix)\in \vU_4$. This will
reduce the number of cases to consider significantly.

\bigskip
\noindent \textbf{I}. $a>0,b,c,d\ge0,e=0$.
\begin{enumerate}[{1)}]
\item\label{case:a>0b=0c=0d=0e=0} $a>0,b=c=d=e=0\Rightarrow \ga=0, \gb,\gc,\gd,\gg\ge0$. Then
\begin{equation*}
  1-C (x-1)^\gb (x+1)^\gc (x-i)^\gd (x+i)^\gg=D x^a.
\end{equation*}
Clearly $C=-D$ and $a=\gb+\gc+\gd+\gg$.

i) If $a=1$, then $\gb+\gc+\gd+\gg=1$.

\begin{itemize}
 \item $\gb=1 \Rightarrow \gc=\gd=\gg=0$. Then $x=1\Rightarrow D=1 \Rightarrow \boxed{f=x}.$

 \item $\gc=1 \Rightarrow \gb=\gd=\gg=0$. Then $x=-1\Rightarrow D=-1 \Rightarrow\boxed{f=-x}.$

 \item $\gd=1 \Rightarrow \gb=\gc=\gg=0$. Then $x=i\Rightarrow D=-i \Rightarrow\boxed{f=-ix}.$

 \item $\gg=1 \Rightarrow \gb=\gc=\gd=0$. Then $x=-i\Rightarrow D=i \Rightarrow\boxed{f=ix}.$
\end{itemize}

ii) If $a=\gb+\gc+\gd+\gg\ge 2\Rightarrow \gb-\gc-i(\gd-\gg)=0$ by 2HD. So $\gb=\gc,\gd=\gg$ and the equation becomes $1-C(x^2-1)^\gb(x^2+1)^\gc=Dx^a$ which yields three subcases:
 \begin{itemize}
   \item  $\gb=0,\gc=1$. Then  $\boxed{f=-x^2}.$

   \item $\gb=1,\gc=0$. Then $\boxed{f=x^2}.$

   \item  $\gb+\gc\ge2$. Then $\gb=\gc=1$ by 2HD. Thus $\boxed{f=x^4}.$
 \end{itemize}

\item \label{case:a>0b>0c=0d=0e=0} $a>0,b>0,c=d=e=0\Rightarrow \ga=\gb=0, \gc\ge 0, \gd\ge 0, \gg\ge0$. Then
\begin{equation*}
  1-C (x+1)^\gc(x-i)^\gd (x+i)^\gg=D x^a (x-1)^b.
\end{equation*}
This is impossible by the 2-adic argument if we consider $x=0$ and $x=1$.

\item \label{case:a>0b=0c>0d=0e=0} $a>0,b=0,c>0,d=e=0\Rightarrow \ga=\gc=0, \gb\ge0, \gd\ge0, \gg\ge0$. Then
\begin{equation*}
  1-C (x-1)^\gb(x-i)^\gd (x+i)^\gg=D x^a (x+1)^c.
\end{equation*}
This is impossible by the 2-adic argument by taking $x=0$ and $x=-1$.

\item \label{case:a>0b=0c=0d>0e=0} $a>0,b=c=0,d>0,e=0\Rightarrow \ga=\gd=0, \gb\ge0, \gc\ge0, \gg\ge0$. Then
\begin{equation*}
  1-C (x-1)^\gb (x+1)^\gc(x+i)^\gg=D x^a (x-i)^d.
\end{equation*}
This is impossible by the 2-adic argument by taking $x=0$ and $x=i$.

\item\label{case:a>0b>0c>0d=0e=0} $a>0,b>0,c>0,d=e=0 \Rightarrow \ga=\gb=\gc=0, \gd\ge 0, \gg\ge0$.
\begin{equation*}
  1-C (x-i)^\gd (x+i)^\gg=D x^a (x-1)^b (x+1)^c.
\end{equation*}
This is impossible by the 2-adic argument by taking $x=0$ and $x=1$.

\item \label{case:a>0b>0c=0d>0e=0} $a>0,b>0,c=0,d>0,e=0\Rightarrow \ga=\gb=\gd=0, \gc\ge 0, \gg\ge0$.
\begin{equation*}
  1-C(x+1)^\gc (x+i)^gg=D x^a (x-1)^b (x-i)^d .
\end{equation*}
This is impossible by the 2-adic argument by taking $x=0$ and $x=i$.

\item \label{case:a>0b=0c>0d>0e=0} $a>0,b=0,c>0,d>0,e=0 \Rightarrow \ga=\gc=\gd=0, \gb\ge 0, \gg\ge0$.
\begin{equation*}
  1-C (x-1)^\gb  (x+i)^\gg=D x^a (x+1)^c(x-i)^d.
\end{equation*}
This is impossible by the 2-adic argument by taking $x=0$ and $x=i$.

\item \label{case:a>0b>0c>0d>0e=0} $a>0,b>0,c>0,d>0,e=0 \Rightarrow \ga=\gc=\gb=0, \gg\ge0$.
\begin{equation*}
  1-C (x+i)^\gg=D x^a(x-1)^b(x+1)^c(x-i)^d.
\end{equation*}
This is impossible by the 2-adic argument by taking $x=0$ and $x=1$.

\newpage
\noindent\ \hskip-1.1cm \textbf{II}. $a>0,b<0,e=0$.

\item \label{case:a>0b<0c=0d=0e=0} $a>0,b<0,c=0,d=e=0\Rightarrow \ga=0, \gb=b=-k, \gc\ge0, \gd\ge0, \gg\ge0$. Then
\begin{equation*}
 (x-1)^k-C (x+1)^\gc(x-i)^\gd  (x+i)^\gg=D x^a .
\end{equation*}

First $x=0 \Rightarrow C=\pm 1,\pm i$.
Note that $\gc\ne0, x=-1\Rightarrow \ord_2(D)=k$; $\gd\ne0, x=i\Rightarrow \ord_2(D)=k/2$; $\gg\ne0, x=-i\Rightarrow \ord_2(D)=k/2$.
So we have four cases to consider.

i) $\gd=\gg=0.$ First $x=0\Rightarrow C=(-1)^k$. Rewriting it as $(1-x)^k=(1+x)^\gc+(-1)^k D x^a$
we see that $(-1)^k D<0$ and $a=1$ since
the linear term is $-kx$ on the LHS. Moreover, if $k> 2$ then there are at least two negative terms on the left
which is impossible. If $k=1$ then $\gc\le 1$ while if $k=2$ then $\gc=2$.

\begin{itemize}
  \item If  $k=1$, $\gc=0$ then $D=1 \Rightarrow \boxed{f=\frac{x}{x-1}}.$

  \item If  $k=1$, $\gc=1$ then $D=2\Rightarrow \boxed{f=\frac{2x}{x-1}}.$

  \item If  $k=2$ then $\gc+D=-2\Rightarrow D=-4 \Rightarrow \boxed{f=\frac{-4x}{(x-1)^2}}.$
\end{itemize}

ii) $\gd=0,\gg\ge 1, \gc=0$.
\begin{itemize}
  \item $k\ge a+1\ge 2$. Then $k=\gg$ and $C=1$. If $a<k-1$ then by 2HD $k-\gg i=0$ which is absurd. So $a=k-1$ and $k-l i=D$.
  Since $k\ge 2$ then 3HD $\binom{k}{2}-\binom{k}{2}i^2=0$ which is impossible.

  \item $k=a$. If $a=\gg+1\ge 2$ then $D=1$ and by 2HD $k=C \Rightarrow k=1$ which is absurd. $a\ge \gg+2\ge 3$
  then $D=1$ and by 2HD $k=0$, again impossible. So we must have $a=\gg$ and by 2HD $k=C k i \Rightarrow C=-i, D=1+i.$
  Note that  $k\ge 2$ is impossible since otherwise by 3HD we would get $\binom{k}{2}=-i\binom{k}{2}i^2$. Hence $\boxed{f=\frac{(1+i)x}{x-1}}.$
\end{itemize}

iii) $\gd\ge 1,\gg=0, \gc=0$. Taking complex conjugation of the coefficients in case ii) above we get $\boxed{f=\frac{(1-i)x}{x-1}}.$

iv) $\gd\ge 1,\gg\ge 1, \gc=0$. If $\gd+\gg>k$ then $a=k$ and $C=-D$ which contradicts to $\ord_2(D)=k/2$. Thus $\gd+\gg=k$.
\begin{itemize}
  \item If $a=k$ then $1-C=D$ and by 2HD $k+C(\gd-\gg)i=0$ which is impossible for $C=\pm1, \pm i$.

  \item If $a=k-1$ then $C=1$. If $a\ge 2$ then from the constant term and the linear term we say that
  $1=(-i)^\gd i^\gg$ and $k=\gd(-i)^{\gd-1}i^\gg+\gg (-i)^\gd i^{\gg-1}=(\gd-\gg)i$ which is absurd. So $a=1, k=2, \gd=\gg=1\Rightarrow\boxed{f=\frac{-2x}{(x-1)^2}}$.

  \item  If $a<k-1$ then $C=1$ and $k+(\gd-\gg)i=0$ by 2HD, which is absurd.
\end{itemize}

\item \label{case:a>0b<0c=0d>0e=0} $a>0,b<0,c=0,d>0,e=0\Rightarrow \ga=0, \gb=b=-k, \gc\ge0,\gd=0,\gg\ge0$. Then
\begin{equation*}
 (x-1)^k-C (x+1)^\gc  (x+i)^\gg=D x^a (x-i)^d.
\end{equation*}
First $x=0\Rightarrow C=\pm 1, \pm i$; $x=i\Rightarrow k=\gc+2\gg$ by the 2-adic argument.
Similarly, $\gc\ne 0, x=-1 \Rightarrow \ord_2(D)=k-d/2$ while $\gg\ne 0, x=-i \Rightarrow \ord_2(D)=k/2-d$.
Thus $\gc\gg=0$.

i)  $\gg=0$. Then $k=\gc$.

\begin{itemize}
  \item $a+d=k, 1-C=D, k+C\gc=D d i=(C-1)d i$ by 2HD. Then $C=-i$ is the only possible choice.
So $d=\gc=k>0$ which contradicts $k=a+d$.

  \item $a+d=k-1\ge 2$ then  $C=1$, $2k=-D$ by 2HD and $0=-Dd i=2kd i$ by 3HD, which is impossible.

  \item $a+d<k-1$ then $k+\gc=0$ by 2HD which is absurd.
\end{itemize}

ii) $\gc=0$. Then $k=2\gg>\gg$ which implies that $D=1$ and $a+d=k$. By 2HD $k=di$ if $\gg<k-1$, which is absurd.
Thus $\gg=k-1$ and $k+C=di$ which cannot hold for $C=\pm 1, \pm i$.

\item\label{case:a>0b<0c>0d=0e=0} $a>0,b<0,c>0,d=e=0\Rightarrow \ga=\gc=0, \gb=b=-k, \gd\ge0, \gg\ge0$. Then
\begin{equation*}
 (x-1)^k-C (x-i)^\gd  (x+i)^\gg=D x^a (x+1)^c.
\end{equation*}

Note that $x=0 \Rightarrow C=\pm1,\pm i$ and $x=-1\Rightarrow 2k=\gd+\gg$ by the 2-adic argument.
Thus $a+c=\gd+\gg \ge 2$ and $C=-D$. Further $\gd+\gg\ge k+1$.
If $\gd+\gg\ge k+2$ then $(\gd-\gg)i=c$ by 2HD, which is absurd.
So $2k=\gd+\gg\ge k+1\Rightarrow k=a=c=1, \gd+\gg=2 \Rightarrow -k-C(\gg-\gd)i=Dc=-Cc$ by 2HD, which cannot hold
for $C=\pm1,\pm i$ if $\gg\ne \gd$. Thus $\gg=\gd=1$ and the equation becomes
$ x-1+D(x^2+1)=Dx(x+1)\Rightarrow D=-1, \boxed{f=\frac{x(x+1)}{x-1}}.$

\item\label{case:a>0b<0c>0d>0e=0} $a>0,b<0,c>0,d>0,e=0\Rightarrow \ga=0, \gb=b=-k, \gc=0, \gd=0,\gg\ge0$. Then
\begin{equation*}
 (x-1)^k-C (x+i)^\gg=D x^a (x+1)^c (x-i)^d.
\end{equation*}
This is impossible by the 2-adic argument. Indeed, $x=0\Rightarrow C=\pm 1, \pm i$;  $x=-1\Rightarrow 2k=\gg$
while $x=i\Rightarrow k=2\gg$ which is absurd.

\item \label{case:a>0b<0c<0d=0e=0} $a>0,b<0,c<0,d=e=0\Rightarrow \ga=0, \gb=b=-k, \gc=c=-l, \gd\ge0,\gg\ge0$. Then
\begin{equation*}
 (x-1)^k (x+1)^l-C (x-i)^\gd (x+i)^\gg=D x^a.
\end{equation*}

i) $\gd+\gg=0$. Then $k=l=1$ by the sign pattern and we get $\boxed{f=\frac{x^2}{x^2-1}}.$

ii) $\gd+\gg\ge 1$. Then $x=0\Rightarrow C=\pm 1, \pm i$; $x=1\Rightarrow \ord_2(D)=(\gd+\gg)/2$;
$x=i$ or $x=-i\Rightarrow \ord_2(D)=(k+l)/2$. Thus $k+l=\gd+\gg\ge 2$.

\begin{itemize}
\item $a=1,k+l=\gd+\gg\ge 3\Rightarrow C=1$. Then $l-k+(\gd-\gg)i=0$ by 2HD. So $(x^2-1)^k-(x^2+1)^k=Dx$
 which is absurd since the LHS is even function while the RHS is odd.

\item $a=1,k+l=\gd+\gg\Rightarrow k=l=1$. Then $x^2-1-(x-i)^\gd (x+i)^\gg=Dx$. Thus $\gd=2,\gg=0$ or $\gd=0,\gg=2$. So we get
$\boxed{f=\frac{2ix}{x^2-1}}$ or $\boxed{f=\frac{-2ix}{x^2-1}}$.

\item $a\ge 2, k+l=\gd+\gg=a \Rightarrow 1-C=D$. Then $l-k+(\gd-\gg)i=0$ by 2HD $\Rightarrow(x^2-1)^k-C(x^2+1)^k=Dx^{2k}$.
The constant term is $(-1)^k-C=0$ so we must have $k$ is odd and $C=-1$. By the sign pattern we see that $k=1\Rightarrow
\boxed{f=\frac{2x^2}{x^2-1}}$.

\item $a\ge 2, k+l=\gd+\gg>a \Rightarrow C=1$. The equation becomes
$(x-1)^k (x+1)^l-(x-i)^\gd (x+i)^\gg=D x^a.$
The constant term $(-1)^k- (-i)^\gd i^\gg=0$ and the linear term
$l-(-1)^k k-(-1)^k((-i)^{-1}\gd+i^{-1}\gg)=l-(-1)^k k-(-1)^k(\gd-\gg)i=0$. Thus $\gd=\gg=l=k$ and $k$ is even.
Thus $(x^2-1)^k-(x^2+1)^k=D x^a.$ By the sign pattern we must have $k=2,a=2$ and $D=-4\Rightarrow \boxed{f=\frac{-4x^2}{(x^2-1)^2}}.$
\end{itemize}

\item\label{case:a>0b<0c=0d<0e=0} $a>0,b<0,c=0,d<0,e=0\Rightarrow \ga=0, \gb=b=-k, \gc\ge0, \gd=d=-m,\gg\ge0$. Then
\begin{equation*}
 (x-1)^k(x-i)^m -C(x+1)^\gc (x+i)^\gg=D x^a.
\end{equation*}
Note $x=0\Rightarrow C=\pm1,\pm i$; $x=1\Rightarrow \ord_2(D)=\gc+\gg/2$;
 $x=i \Rightarrow \ord_2(D)=\gc/2+\gg$. So we get $\gc=\gg$. By 2HD $-k-mi=C\gc(1+i)\Rightarrow C=-1\Rightarrow k=m=\gc=\gg$. Thus $(x^2-(i+1)x+i)^k-(x^2+(i+1)x+i)^k=D x^a$. So $k=1$ and $D=-2(i+1)\Rightarrow \boxed{f=\frac{-2(1+i)x}{(x-1)(x-i)}}.$

\item\label{case:a>0b<0c<0d>0e=0} $a>0,b<0,c<0,d>0,e=0\Rightarrow \ga=0, \gb=b=-k, \gc=c=-l, \gd=0,\gg\ge0$. Then
\begin{equation*}
 (x-1)^k (x+1)^l-C  (x+i)^\gg=D x^a(x-i)^d.
\end{equation*}
Note that $x=0\Rightarrow C=\pm1,\pm i$ by the 2-adic argument.

i)  $k+l<\gg-1$. Then $\gg=a+d, C=-D$. By 2HD, $-C\gg i=-C(a+d) i=-d D i=C d i$ which is absurd.

i)  $k+l=\gg-1$. Then $\gg=a+d, C=-D$. By 2HD, $l-k-C\gg i=-C(a+d) i=-d D i=C d i\Rightarrow l-k=C(a+2d)i$ which cannot hold for $C=\pm1,\pm i$.

iii) $k+l=\gg\ge 2$. Then $x=1 \Rightarrow \ord_2(D)=\gg/2-d/2$. If $\gg>0$ then $x=-i \Rightarrow \ord_2(D)=(k+l)/2-d=\gg/2-d$ which is absurd.

iv) $k+l=\gg+1\ge 2$. Thus $D=1,k+l=a+d$. By 2HD, $l-k+C=-d i$ which can hold only if $d=1, C=-i, l=k=1\Rightarrow a=d=1 \Rightarrow \boxed{f=\frac{x(x-i)}{x^2-1}}$.

v) $k+l>\gg+1\ge 2$. Thus $D=1,k+l=a+d$. By 2HD, $l-k=-d i$ which is absurd.

\item \label{case:a>0b<0c>0d<0e=0} $a>0,b<0,c>0,d<0,e=0\Rightarrow \ga=\gc=0, \gb=b=-k, \gd=d=m, \gg\ge0$. Then
\begin{equation*}
   (x-1)^k (x-i)^m -C (x+i)^\gg=D x^a (x+1)^c.
\end{equation*}
Note that $2k+m=\gg$ and $C=\pm1,\pm i$ by the 2-adic argument.
Now $x=1 \Rightarrow \ord_2(D)=\gg/2-c$ while $x=i \Rightarrow \ord_2(D)=\gg-c/2$ which is absurd.

\item\label{case:a>0b<0c<0d<0e=0} $a>0,b<0,c<0,d<0,e=0\Rightarrow \ga=0, \gb=b=-k, \gc=c=-l, \gd=d=-m,\gg\ge0$. Then
\begin{equation*}
 (x-1)^k (x+1)^l(x-i)^m-C  (x+i)^\gg=D x^a.
\end{equation*}
First $x=0 \Rightarrow C=\pm1,\pm i$. Now $x=1 \Rightarrow \ord_2(D)=\gg/2$ while $x=i \Rightarrow \ord_2(D)=\gg$. Thus $\gg=0$.
So $a=k+l+m$. By 2HD we get $-k+l-mi=0$ which is absurd.

\newpage
\noindent\ \hskip-1.1cm \textbf{III}. $a>0,b\ge0,e=0,c<0$ or $d<0$.

\item \label{case:a>0b>0c<0d<0e=0} $a>0,b>0,c<0,d<0,e=0\Rightarrow \ga=\gb=0, \gc=c=-l, \gd=d=-m,\gg\ge0$. Then
\begin{equation*}
 (x+1)^l(x-i)^m-C (x+i)^\gg=D x^a(x-1)^b.
\end{equation*}
Taking $x\leftrightarrow -x$ and then applying complex conjugation on the coefficients we see
that there is no solution by Case \ref{case:a>0b<0c>0d<0e=0}.

\item \label{case:a>0b>0c<0d=0e=0} $a>0,b>0,c<0,d=e=0\Rightarrow \ga=0, \gb=0, \gc=c=-l, \gd\ge 0,\gg\ge0$. Then
\begin{equation*}
 (x+1)^l-C(x-i)^\gd (x+i)^\gg=D x^a(x-1)^b.
\end{equation*}
By Case \ref{case:a>0b<0c>0d=0e=0}) using substitution $x\leftrightarrow -x$ we get $\boxed{f=\frac{x(1-x)}{x+1}}.$

\item \label{case:a>0b>0c=0d<0e=0} $a>0,b>0,c=0,d<0,e=0\Rightarrow \ga=\gb=0,\gc\ge 0,\gd=d=-m,\gg\ge0$. Then
\begin{equation*}
 (x-i)^m-C (x+1)^\gc (x+i)^\gg=D x^a(x-1)^b.
\end{equation*}
Applying $x\leftrightarrow ix$ and then taking complex conjugation on the coefficients we see
that there is no solution by Case \ref{case:a>0b<0c=0d>0e=0}.

\item \label{case:a>0b>0c<0d>0e=0} $a>0,b>0,c<0,d>0,e=0\Rightarrow \ga=0, \gb=0, \gc=c=-l, \gd=0,\gg\ge0$. Then
\begin{equation*}
 (x+1)^l-C(x+i)^\gg=D x^a(x-1)^b (x-i)^d.
\end{equation*}
Applying $x\leftrightarrow -x$ and then taking complex conjugation on the coefficients we see
that there is no solution by Case \ref{case:a>0b<0c>0d>0e=0}).

\item \label{case:a>0b>0c>0d<0e=0} $a>0,b>0,c>0,d<0,e=0\Rightarrow \ga=\gb=\gc=0, \gd=d=-m,\gg\ge0$. Then
\begin{equation*}
(x-i)^m-C(x+i)^\gg=D x^a(x-1)^b (x+1)^c .
\end{equation*}
Note that $m=\gg$, $|C|=1$ by 2-adic argument.
Further $a=1$ since otherwise the constant term and
linear terms are both 0 on the LHS which is impossible. Also, $C=1$ or $-1$ depending whether $m$ is even or odd.

i) $m$ even, $C=1$ so $1+b+c\le m-1$.

\begin{itemize}
  \item $1+b+c< m-1$. Then by 2HD $-2m i=0$ which is absurd.

  \item $1+b+c=m-1$. By 2HD is $-2m i=D$ so that all coefficients on the RHS are imaginary. By 3HD coefficient $2\binom{m}{3} i=-2mi\big(\binom{b}{2}-\binom{c}{2}\big) \Rightarrow b=c$. The equation becomes
$(x-i)^{2q}-(x+i)^{2q}=-4qi x(x^2-1)^{q-1}$. By 4HD, $2\binom{2q}{3} i=4qi(q-1)\Rightarrow q=2 \Rightarrow \boxed{f=\frac{8i x(1-x^2)}{(x-i)^4}}.$
\end{itemize}

ii) $m$ is odd, $C=-1$. Then $1+b+c=m$ and $D=2$. By 2HD we see that $-2m i=2(c-b) \Rightarrow b=c, m=0$ which is absurd.

\item \label{case:a>0b=0c<0d<0e=0} $a>0,b=0,c<0,d<0,e=0\Rightarrow \ga=0, \gb\ge0, \gc=c=-l, \gd=d=-m,\gg\ge0$. Then
\begin{equation*}
 (x+1)^l (x-i)^m-C (x-1)^\gb (x+i)^\gg=D x^a.
\end{equation*}
Applying $x\leftrightarrow -x$ in Case \ref{case:a>0b<0c=0d<0e=0}) and then taking complex conjugation on the coefficients we get $\boxed{f=\frac{2(1-i)x}{(x+1)(x-i)}}.$

\item \label{case:a>0b=0c<0d=0e=0} $a>0,b=0,c<0,d=0,e=0\Rightarrow \ga=0, \gb\ge0, \gc=c=-l, \gd\ge 0,\gg\ge0$. Then
\begin{equation*}
 (x+1)^l-C (x-1)^\gb (x-i)^\gd (x+i)^\gg=D x^a.
\end{equation*}
By Case \ref{case:a>0b<0c=0d=0e=0}) we get
\begin{equation*}
\boxed{f=\frac{x}{x+1},\quad  \frac{2x}{x+1},\quad \frac{4x}{(x+1)^2},\quad \frac{(1+i)x}{x+1},\quad\frac{(1-i)x}{x+1},\quad \frac{2x}{(x+1)^2}.}
\end{equation*}

\item \label{case:a>0b=0c=0d<0e=0} $a>0,b=0,c=0,d<0,e=0\Rightarrow \ga=0,\gb\ge0,\gc\ge0, \gd=d=-m,\gg\ge0$. Then
\begin{equation*}
 (x-i)^m-C (x-1)^\gb(x+1)^\gc(x+i)^\gg=D x^a .
\end{equation*}
Applying $x \to ix$ in Case \ref{case:a>0b=0c<0d=0e=0}) we get
\begin{equation*}
\boxed{f=\frac{x}{x-i},\quad  \frac{2x}{x-i},\quad \frac{-4ix}{(x-i)^2},\quad \frac{(1+i)x}{x-i},\quad\frac{(1-i)x}{x-i},\quad \frac{-2ix}{(x-i)^2}.}
\end{equation*}

\item \label{case:a>0b=0c<0d>0e=0} $a>0,b=0,c<0,d>0,e=0\Rightarrow \ga=0, \gb\ge0, \gc=c=-l, \gd=0,\gg\ge0$. Then
\begin{equation*}
 (x+1)^l -C (x-1)^\gb  (x+i)^\gg=D x^a (x-i)^d.
\end{equation*}
By Case \ref{case:a>0b<0c=0d>0e=0}) there is no solution by $x\leftrightarrow -x$ and applying complex conjugation on coefficients.

\item \label{case:a>0b=0c>0d<0e=0} $a>0,b=0,c>0,d<0,e=0\Rightarrow \ga=0, \gb\ge0, \gc=0, \gd=d=-m,\gg\ge0$. Then
\begin{equation*}
 (x-i)^m -C (x-1)^\gb  (x+i)^\gg=D x^a (x+1)^c.
\end{equation*}
By Case \ref{case:a>0b>0c=0d<0e=0}) there is no solution by $x\leftrightarrow -x$ and applying complex conjugation on coefficients.

\bigskip
\noindent\ \hskip-1.1cm \textbf{IV}. $a=0,b,c,d\ge0,e=0$. We set $1/D=-B$. Then we must have $C=-D$.

\item \label{case:a=0b=0c=0d=0e=0} $a=0,b=c=d=e=0\Rightarrow \ga,\gb,\gc,\gd,\gg\ge0$. This is clearly impossible.

\item \label{case:a=0b>0c=0d=0e=0} $a=0,b>0,c=d=e=0\Rightarrow \ga\ge 0,\gb=0, \gc\ge 0, \gd\ge 0,\gg\ge0$.
Then we get
\begin{equation*}
 B-x^\ga (x+1)^\gc (x-i)^\gd (x+i)^\gg=-(x-1)^b.
\end{equation*}

i) $\ga\ge 1$. $x=0\Leftrightarrow |B|=1$ which is impossible by the 2-adic argument except for $\gc=\gd=\gg=0$
which yield $b=1$ by the sign pattern. Thus $\ga=B=1\Rightarrow \boxed{f=1-x}.$

ii) $\ga=0$. $\gc\ne 0,x=-1\Rightarrow \ord_2(B)=b$; $x=i\Rightarrow \ord_2(B)=b/2$ if $\gd\ne 0$;
$x=-i\Rightarrow \ord_2(B)=b/2$ if $\gg\ne 0$. Thus $\gc(\gd+\gg)=0$.
\begin{itemize}
 \item $\gc=1,\gd=\gg=0.$ Then $b=1$ since the linear term on the LHS has positive coefficient. Thus we get $\boxed{f=\frac{1-x}{2}}.$

 \item $\gc=0$. If $\gd+\gg\ge 2$ then by 2HD $(\gg-\gd)i=-b$ which is absurd. Thus $\gd=1,\gg=0\Rightarrow   \boxed{f=\frac{1-x}{1-i}}$ or $\gd=0,\gg=1\Rightarrow \boxed{f=\frac{1-x}{1+i}}$.
\end{itemize}

\item \label{case:a=0b=0c>0d=0e=0} $a=0,b=0,c>0, d=e=0\Rightarrow \ga=\gc=0, \gb\ge0, \gd\ge0,\gg\ge0$. Then we must have
\begin{equation*}
 1-Cx^\ga (x-1)^\gb(x-i)^\gd(x+i)^\gg=D(x+1)^c.
\end{equation*}
By Case \ref{case:a=0b>0c=0d=0e=0}),  applying $x\leftrightarrow -x$ and complex conjugation on coefficients we get
\begin{equation*}
\boxed{f=x+1, \quad\frac{x+1}{2}, \quad \frac{x+1}{1+i},  \quad \frac{x+1}{1-i}}.
\end{equation*}

\item \label{case:a=0b=0c=0d>0e=0} $a=0,b=c=0,d>0,e=0\Rightarrow \ga\ge0, \gb\ge0,\gc\ge0, \gd=0,\gg\ge0$. By symmetry $i\leftrightarrow -i$ from above we get
\begin{equation*}
  1-Cx^\ga (x-1)^\gb (x+1)^\gc (x+i)^\gg=D(x-i)^d.
\end{equation*}
Applying $x\leftrightarrow -ix$ in Case \ref{case:a=0b>0c=0d=0e=0}) we get
\begin{equation*}
\boxed{f=\frac{x-i}{-i}, \quad\frac{x-i}{-2i}, \quad \frac{i-x}{i+1}, \quad \frac{x-i}{1-i}}.
\end{equation*}

\item \label{case:a=0b>0c>0d=0e=0} $a=0,b>0,c>0, d=e=0 \Rightarrow \ga=\gb=\gc=0, \gd\ge 0,\gg\ge0$. Then we must have
\begin{equation*}
  B-x^\ga (x-i)^\gd (x+i)^\gg=-(x-1)^b (x+1)^c.
\end{equation*}
Note that $x=0\Rightarrow B=\pm 1$ if $\ga\ne0$; $x=i\Rightarrow \ord_2(B)=(b+c)/2$ if $\gd\ne 0$;
$x=-i\Rightarrow \ord_2(B)=(b+c)/2$ if $\gg\ne 0$. Thus $\ga(\gd+\gg)=0$.

i) $\ga\ge 2, \gd=\gg=0$. So $x=0\Rightarrow B=(-1)^b$. We also have $b=c$ since the linear term vanishes on the LHS. Then $(-1)^b-x^{2b}=-(x^2-1)^b$. Thus $b=1$ since the LHS has only two terms $\Rightarrow \boxed{f=1-x^2}.$

ii) $\ga=0, \gd+\gg\ge 1$. If $\gd+\gg\ge 2$ then by 2HD $(\gd-\gg)i=b-c\Rightarrow \gd=\gg, b=c$. So $B-(x^2+1)^\gd=-(x^2-1)^b$.
by the sign pattern we get $\gb=b=1\Rightarrow \boxed{f=\frac{1-x^2}2}.$

\item \label{case:a=0b>0c=0d>0e=0} $a=0,b>0,c=0,d>0,e=0\Rightarrow \ga\ge0, \gb=0, \gc\ge0, \gd=0,\gg\ge0$. Then we must have
\begin{equation*}
  B-x^\ga (x+1)^\gc (x+i)^\gg=-(x-1)^b (x-i)^d.
\end{equation*}
Note that $x=1\Rightarrow \ord_2(B)=\gc+\gg/2$ while $x=i\Rightarrow \ord_2(B)=\gc/2+\gg$
which implies that $\gc=\gg=0$. Thus $\ga\ne0$ so that $x=0\Rightarrow |B|=1$;
But there are at least three nonzero terms on the RHS, constant, linear and the highest degree term.
This contradicts to the fact that LHS has only two terms. Hence, this case has no solution.

\item \label{case:a=0b=0c>0d>0e=0} $a=0,b=0,c>0, d>0,e=0 \Rightarrow \ga\ge 0, \gb\ge 0, \gc=\gd=0, \gg\ge0$. Then we must have
\begin{equation*}
  1-Cx^\ga (x-1)^\gb (x+i)^\gg =D(x+1)^c(x-i)^d.
\end{equation*}
This has no solution by setting $x\Leftrightarrow -x$ and applying complex conjugation in Case \ref{case:a=0b>0c=0d>0e=0}).

\item \label{case:a=0b>0c>0d>0e=0} $a=0,b>0,c>0,d>0,e=0 \Rightarrow \ga=\gc=\gb=0, \gg\ge0$. Then we must have
\begin{equation*}
  B- x^\ga(x+i)^\gg =-(x-1)^b (x+1)^c (x-i)^d.
\end{equation*}
By 2-adic argument we see that $\gg=0$ and $\ga=b+c+d\ge3$. But this is absurd since the RHS has at least three nonzero terms:
constant, linear and the highest degree term..

\newpage
\noindent\ \hskip-1.1cm \textbf{V}. $a=0,b<0,e=0$.

\item \label{case:a=0b<0c=0d=0e=0} $a=0,b<0,c=0,d=e=0\Rightarrow \ga\ge 0, \gb=b=-k, \gc\ge0, \gd\ge0,\gg\ge0$. Then
\begin{equation*}
 (x-1)^k-C x^\ga (x+1)^\gc (x-i)^\gd  (x+i)^\gg=D .
\end{equation*}
By Case \ref{case:a=0b>0c=0d=0e=0}) we get
\begin{equation*}
 \boxed{f=\frac{1}{1-x},\qquad \frac{2}{1-x},\qquad  \frac{1+i}{1-x},  \qquad \frac{1-i}{1-x}}.
\end{equation*}

\item \label{case:a=0b<0c=0d>0e=0} $a=0,b<0,c=0,d>0,e=0\Rightarrow \ga\ge 0, \gb=b=-k, \gc\ge0,\gd=0,\gg\ge0$. Then
\begin{equation*}
 (x-1)^k-C x^\ga (x+1)^\gc  (x+i)^\gg=D (x-i)^d.
\end{equation*}
Note that $x=-1\Rightarrow \ord_2(D)=k-d/2$ if $\gc\ne0$; $x=-i\Rightarrow \ord_2(C)=k/2-d$ if $\gg\ne0$. Thus $\gc\gg=0$.

i) $\ga=\gc=\gg=0$. Then by 2HD $k=Dd i$ so $d=1$ since otherwise by 3HD $\binom{k}{2}=-D\binom{d}{2}$ which is absurd.
Thus $d=1\Rightarrow\boxed{f=\frac{x-i}{x-1}}.$

ii) If $\ga\ge 1$, $\gc=\gg=0$ then by Case \ref{case:a>0b<0c=0d=0e=0}) we find one solution $\boxed{f=\frac{1+ix}{1-x}}.$

iii) If $\ga=0,\gc\ge1,\gg=0$ then $ (x-1)^k-C (x+1)^\gc=D (x-i)^d$. By comparing the constant term we see that $C=(-1)^k-D(-i)^d$.
\begin{itemize}
  \item $k>\gc+1$. Then $D=1$ and $d=k$. By 2HD $-k=-di$ which is impossible.

  \item $k=\gc+1\ge 2$. Then $D=1$ and $d=k$. By 2HD $k i-k=C=(-1)^k(1-i^k)$ which has no solution.

  \item $k=\gc$. If $d=k$ then $1-C=D$ and, if $k\ge2$, $1-C=-D$ by 3HD which is absurd. Thus $k=1$,
   $1-C=-D i\Rightarrow C=i, D=1-i\Rightarrow \boxed{f=\frac{(1-i)(x-i)}{x-1}}.$

  \item $k=\gc-1$. Then $\gc=d$, $C=-D$ The constant term yields $C=(-1)^k+C(-i)^k\Rightarrow C=1/2, \frac{-1}{1\pm i}$,
   none of which satisfies $1-Cd=-Ddi=Cdi$ by 2HD.

  \item $k<\gc-1$. Then $\gc=d$, $C=-D$ and by 2HD $-Cd=-Ddi$ which is impossible.
\end{itemize}

iv) If $\ga=0,\gc=0,\gg\ge1$ then similarly to Case v) above, we have $\gg-1\le k\le \gg+1$. By comparing the constant term we see that $(-1)^\gg C=(-1)^k-D(-i)^d$.
\begin{itemize}
  \item $k=\gg+1\ge2$. Then $D=1$ and $d=k$. By 2HD $k i-k=C$, and by constant term $(-1)^k-C i^{k-1}=(-i)^k$ which has no solution.

  \item $k=\gg\ge2$. If $d=k$ then $1-C=D$ and by 2HD $k+Cki=D k i \Rightarrow 1+Ci=(1-C)i \Rightarrow C=\frac{1+i}2,D=\frac{1-i}2$.
  If $k\ge2$ then $1+C=-D$ by 3HD which is absurd. Thus $k=1 \Rightarrow \boxed{f=\frac{(1-i)(x-i)}{2(x-1)}}.$

  \item $k=\gg-1$. Then $\gg=d$, $C=-D$ The constant term yields $C=-1+Ci^k\Rightarrow C=-1/2, \frac{-1}{1\pm i}$,
   none of which satisfies $1-Cdi=-Ddi=Cdi$ by 2HD.
\end{itemize}

v) If $\ga\ge 1,\gc=0,\gg\ge1$ then no solution exists by taking complex conjugation in Case \ref{case:a>0b<0c=0d>0e=0}).

vi) If $\ga\ge 1,\gc\ge1,\gg=0$ then there is no solution by Case \ref{case:a>0b<0c>0d=0e=0}).

\item \label{case:a=0b<0c>0d=0e=0} $a=0,b<0,c>0,d=e=0\Rightarrow \ga\ge 0, \gb=b=-k, \gc=0,\gd\ge0,\gg\ge0$. Then
\begin{equation*}
 (x-1)^k-C x^\ga(x-i)^\gd (x+i)^\gg=D (x+1)^c.
\end{equation*}

i) $\ga=\gd=\gg=0$. Then $D=k=c=1$ since otherwise the 2HD cannot match $\Rightarrow \boxed{f=\frac{x+1}{x-1}}.$

ii) If $\ga\ge 1$, $\gd=\gg=0$ then we can rewrite the equation as $(1-x)^k-C' x^\ga =D' (x+1)^c.$
Then $x=0\Rightarrow D'=1$; $x=1\Rightarrow C'=-2^c$ and $\ga=1$ by the sign of the linear term.
\begin{itemize}
  \item If $k=1$ then $c=1\Rightarrow \boxed{f=\frac{1+x}{1-x}}.$
  \item  If $k \ge 2$ then $c=k$ and $-k-C'=k\Rightarrow 2^k=2k \Rightarrow k=2\Rightarrow \boxed{f=\frac{(1+x)^2}{(x-1)^2}}.$
\end{itemize}

iii) If $\ga=0,\gd\ge1,\gg=0$ then $ (x-1)^k-C (x-i)^\gd=D (x+1)^c$. This is essentially the same as Case \ref{case:a=0b<0c=0d>0e=0}.iii) so we get $\boxed{f=\frac{i(x+1)}{x-1}}.$

iv) If $\ga=0,\gd=0,\gg\ge1$ then takeing complex conjugation in v) we get $\boxed{f=\frac{i(x+1)}{1-x}}.$

v) If $\ga\ge 1,\gd=0,\gg\ge1$ then no solution exists by taking complex conjugation in Case \ref{case:a>0b<0c=0d>0e=0}).

vi) If $\ga\ge 1,\gd\ge1,\gg=0$ then there is no solution by Case \ref{case:a>0b<0c=0d>0e=0}).

vii) If $\ga\ge 1,\gd\ge1,\gg\ge 1$ then  $(x-1)^k-C x^\ga (x-i)^\gd (x+i)^\gg=D (x+1)^c.$
\begin{itemize}
  \item $k\ge c+2$. Then $C=1$ and by 2HD $k+(\gg-\gd)i=0$ which is impossible.

  \item $c\ge k+2$. By taking $x\leftrightarrow -x$ this is impossible from above.

  \item $k=c+1$.  Then $C=1,k=\ga+\gd+\gg$ and by 2HD and 3HD $k+(\gg-\gd)i=D, \binom{k}{2}+\binom{\gd}{2}-\gd\gg+\binom{\gg}{2}=Dc$. Thus $D\in\R\Rightarrow D=k, \gg=\gd=k/2$. But by 3HD $\binom{k}{2}+2\binom{\gd}{2}-\gd^2=k(k-1)$ which is impossible.

  \item $c\ge k+1$. By taking $x\leftrightarrow -x$ this is impossible from above.

  \item $k=c$ and $\ga=0$. If $\gd+\gg=k\ge 2$ then $1-C=D$ and by 2HD and 3HD $k+C(\gg-\gd)i=-Dk, \binom{k}{2}+C(\binom{\gd}{2}-\gd\gg+\binom{\gg}{2})=D\binom{k}{2}\Rightarrow 2(\binom{\gd}{2}-\gd\gg+\binom{\gg}{2})=(k-1)(\gg-\gd)i \Rightarrow \gg=\gd$.
      If $\gd=\gg\ge 2$ then this is absurd. So $\gd=\gg=1\Rightarrow k=2, D=-1 \Rightarrow \boxed{f=\frac{-(x+1)^2}{(x-1)^2}}$.

  \item $k=c$ and $\ga=1$. If $\ga+\gd+\gg=k\ge 2$ then by Case \ref{case:a>0b<0c=0d>0e>0}) we get $f=\boxed{f=\frac{(x+1)^4}{(x-1)^4}}$.
\end{itemize}

\item \label{case:a=0b<0c>0d>0e=0} $a=0,b<0,c>0,d>0,e=0\Rightarrow \ga\ge 0, \gb=b=-k, \gc=0, \gd=0,\gg\ge0$. Then
\begin{equation*}
 (x-1)^k-Cx^\ga (x+i)^\gg=D (x+1)^c (x-i)^d.
\end{equation*}

No solution exists by Case \ref{case:a=0b>0c<0d<0e=0}).

\item \label{case:a=0b<0c<0d=0e=0} $a=0,b<0,c<0,d=0,e=0\Rightarrow \ga\ge 0, \gb=b=-k, \gc=c=-l, \gd\ge0,\gg\ge0$. Then
\begin{equation*}
 (x-1)^k (x+1)^l-C x^\ga (x-i)^\gd (x+i)^\gg=D.
\end{equation*}

By Case \ref{case:a=0b>0c>0d=0e=0} we get $\boxed{f=\frac{1}{1-x^2}}$ or $\boxed{f=\frac{2}{1-x^2}}.$

\item\label{case:a=0b<0c=0d<0e=0} $a=0,b<0,c=0,d<0,e=0\Rightarrow \ga\ge 0, \gb=b=-k, \gc\ge0, \gd=d=-m,\gg\ge0$. Then
\begin{equation*}
 (x-1)^k(x-i)^m -C x^\ga (x+1)^\gc (x+i)^\gg=D.
\end{equation*}
Clearly $C=1$. Taking $x=1$ and $x=i$ we see that $\gc=\gg$ by 2-divisibility. Since the RHD is a constant, we get $\gc=\gg\ge 1$. Then taking $x=-1$ and $x=-i$ we see that $k=m$. Thus $\ga=0$ or $\ga=1$. But if $\ga=0$ then
then linear term $-2(k+\gc)(1+i)x$ on the LHS is nonzero which is absurd. If $\ga=1$ then the highest degrees
of the two products on the LHS have different parity which is absurd.

\item \label{case:a=0b<0c<0d>0e=0} $a=0,b<0,c<0,d>0,e=0\Rightarrow \ga\ge 0, \gb=b=-k, \gc=c=-l, \gd=0,\gg\ge0$. Then
\begin{equation*}
 (x-1)^k (x+1)^l-C x^\ga  (x+i)^\gg=D (x-i)^d.
\end{equation*}

i) $\ga=\gg=0$. Then $D=1$ and by 2HD, $-k+l=-di$ which is absurd.

ii) $\ga\ge 1,\gg=0$. Then by Case \ref{case:a>0b<0c<0d=0e=0}) we get $\boxed{f=\frac{(x-i)^2}{x^2-1}}$.

iii) $\ga=0,\gg\ge 1$.
\begin{itemize}
  \item $\gg<d-1$. Then $D=1,k+l=d$ and by 2HD $-k+l=di$ which is absurd.

  \item $\gg=d-1$. Then $D=1,k+l=d$ and by 3HD, $\binom{k}{2}-kl+\binom{l}{2}-C\gg i=-\binom{d}{2}\Rightarrow C\in i\R$. By  2HD $-k+l-C=di\Rightarrow k=l, C=-di$. Thus by 3HD, $-d/2+d(d-1)=-\binom{d}{2}$ which has no solution.

  \item $\gg=d, k+l<d-1$. Then $C=-D$ and by 2HD $C\gg=Dd$ which is absurd.

  \item $\gg=d, k+l=d-1$. Then $C=-D$ and by 2HD $1-Cd i=-Ddi=Cdi\Rightarrow 2 d C i=1\Rightarrow C\in i\R$.
  By 3HD, $-k+l+\binom{d}{2}C=-\binom{d}{2}D\Rightarrow k=l$. By 4HD, $-k-\binom{d}{3}C i^3=\binom{d}{3}D (-i)^3\Rightarrow d=5, k=2, C=-i/10, D=i/10$. But the constant terms do not match.

  \item $\gg=d, k+l=d$. Then $1-C=D$ and by 2HD, $-k+l-Cd i=-Ddi=(C-1)di\Rightarrow 2Cd i=l-k+d i$. By 3HD, $\binom{k}{2}-kl+\binom{l}{2}+\binom{d}{2}C=-\binom{d}{2}D\Rightarrow
      4(\binom{k}{2}-kl+\binom{l}{2})+(d-1)(d+(k-l)i)=(d-1)(-d+(k-l)i)\Rightarrow \binom{k}{2}-kl+\binom{l}{2}+\binom{k+l}{2}=0\Rightarrow k^2+l^2=k+l \Rightarrow k=l=1, d=2, C=D=1/2\Rightarrow
      \boxed{f=\frac{(x-i)^2}{2(x^2-1)}}.$

  \item $\gg=d+1$. Then $C=1,k+l=\gg$. By 2HD and 3HD, $-k+l-\gg i=D, \binom{k}{2}-kl+\binom{l}{2}+\binom{\gg}{2}=-Ddi=(k-l)di-\gg d\Rightarrow k=l, 2k=\gg, -k+\binom{2k}{2}=-2k(2k-1)\Rightarrow k=2/3$ which is absurd.

  \item $\gg>d+1$. Then $C=1,k+l=\gg$. By 2HD and 3HD, $-k+l-\gg i=0$, which is absurd.
\end{itemize}

iv) $\ga\ge 1,\gg\ge 1$. Then taking complex conjugation of coefficients in Case \ref{case:a>0b<0c<0d>0e=0}) we get $\boxed{f=\frac{1+ix}{1-x^2}}$.

\item \label{case:a=0b<0c>0d<0e=0} $a=0,b<0,c>0,d<0,e=0\Rightarrow \ga\ge 0, \gb=b=-k, \gc=0, \gd=d=-m,\gg\ge0$. Then
\begin{equation*}
 (x-1)^k (x-i)^m -C x^\ga (x+i)^\gg=D (x+1)^c.
\end{equation*}
Applying $x \leftrightarrow -x$ and taking complex conjugation of coefficients we see that no solution exists by Case \ref{case:a=0b<0c>0d>0e=0}).

\item \label{case:a=0b<0c<0d<0e=0} $a=0,b<0,c<0,d<0,e=0\Rightarrow \ga\ge 0, \gb=b=-k, \gc=c=-l, \gd=d=-m,\gg\ge0$. Then
\begin{equation*}
 (x-1)^k (x+1)^l(x-i)^m-C x^\ga  (x+i)^\gg=D.
\end{equation*}
No solution exists by Case \ref{case:a=0b>0c>0d>0e=0}).

\bigskip
\noindent\ \hskip-1.1cm \textbf{VI}. $a=0,b\ge 0,c<0,e=0$ or $d<0$.

\item \label{case:a=0b>0c<0d<0e=0} $a=0,b>0,c<0,d<0,e=0\Rightarrow \ga\ge 0, \gb=0, \gc=c=-l, \gd=d=-m,\gg\ge0$. Then
\begin{equation*}
 (x+1)^l(x-i)^m-C x^\ga  (x+i)^\gg=D (x-1)^b.
\end{equation*}
No solution exists by Case \ref{case:a=0b<0c>0d>0e=0}).

\item \label{case:a=0b>0c<0d=0e=0} $a=0,b>0,c<0,d=0,e=0\Rightarrow \ga\ge 0, \gb=0, \gc=c=-l, \gd\ge 0,\gg\ge0$. Then
\begin{equation*}
 (x+1)^l-Cx^\ga (x-i)^\gd  (x+i)^\gg=D (x-1)^b.
\end{equation*}
By Case \ref{case:a=0b<0c>0d=0e=0}) we get
\begin{equation*}
\boxed{f=\frac{x-1}{x+1}, \quad \frac{1-x}{1+x}, \quad \frac{(x-1)^2}{(x+1)^2}, \quad \frac{x-1}{i(x+1)}, \quad \frac{1-x}{i(x+1)}, \quad \frac{-(x-1)^2}{(x+1)^2}, \quad  \frac{(x-1)^4}{(x+1)^4}}.
\end{equation*}

\item \label{case:a=0b>0c=0d<0e=0} $a=0,b>0,c=0,d<0,e=0\Rightarrow \ga\ge 0, \gb=0, \gc\ge 0, \gd=d=-m,\gg\ge0$. Then
\begin{equation*}
 (x-i)^m-Cx^\ga (x+1)^\gc (x+i)^\gg=D (x-1)^b.
\end{equation*}
By Case \ref{case:a=0b<0c=0d>0e=0})  we get $\boxed{f=\frac{1-x}{1+ix},  \quad \frac{x-1}{x-i}, \quad \frac{(1+i)(x-1)}{2(x-i)}, \quad \frac{(1+i)(x-1)}{x-i}}.$

\item \label{case:a=0b>0c<0d>0e=0} $a=0,b>0,c<0,d>0,e=0\Rightarrow \ga\ge 0, \gb=0, \gc=c=-l, \gd=0,\gg\ge0$. Then
\begin{equation*}
 (x+1)^l-C x^\ga  (x+i)^\gg=D (x-1)^b (x-i)^d.
\end{equation*}
There is no solution by Case \ref{case:a=0b<0c>0d<0e=0}).

\item \label{case:a=0b>0c>0d<0e=0} $a=0,b>0,c>0,d<0,e=0\Rightarrow \ga\ge 0, \gb=0, \gc=0, \gd=d=-m,\gg\ge0$. Then
\begin{equation*}
 (x-i)^m-C x^\ga  (x+i)^\gg=D (x-1)^b (x+1)^c.
\end{equation*}
By Case \ref{case:a=0b<0c<0d>0e=0}) we get $\boxed{f=\frac{x^2-1}{(x-i)^2},\  \frac{2(x^2-1)}{(x-i)^2},\ \frac{1-x^2}{1+ix}}.$

\item \label{case:a=0b=0c<0d<0e=0} $a=0,b=0,c<0,d<0,e=0\Rightarrow \ga\ge 0, \gb\ge0, \gc=c=-l, \gd=d=-m,\gg\ge0$. Then
\begin{equation*}
 (x+1)^l (x-i)^m-C x^\ga (x-1)^\gb (x+i)^\gg=D.
\end{equation*}
By Case \ref{case:a=0b=0c>0d>0e=0}) no solution exists.

\item \label{case:a=0b=0c<0d=0e=0} $a=0,b=0,c<0,d=0,e=0\Rightarrow \ga\ge 0, \gb\ge0, \gc=c=-l, \gd\ge 0,\gg\ge0$. Then
\begin{equation*}
 (x+1)^l-C x^\ga (x-1)^\gb (x-i)^\gd (x+i)^\gg=D.
\end{equation*}
By Case \ref{case:a=0b=0c>0d=0e=0}) we get
\begin{equation*}
\boxed{f=\frac{1}{x+1}, \quad\frac{2}{x+1}, \quad \frac{1+i}{x+1}, \quad \frac{1-i}{x+1}}.
\end{equation*}

\item \label{case:a=0b=0c=0d<0e=0} $a=0,b=0,c=0,d<0,e=0\Rightarrow \ga\ge 0, \gb\ge0, \gc\ge 0, \gd=d=-m,\gg\ge0$. Then
\begin{equation*}
 (x-i)^m-C x^\ga (x+1)^\gc (x-1)^\gb  (x+i)^\gg=D.
\end{equation*}
By Case \ref{case:a=0b=0c=0d>0e=0}) we get
\begin{equation*}
\boxed{f=\frac{-i}{x-i}, \quad\frac{-2i}{x-i}, \quad \frac{i+1}{i-x}, \quad \frac{1-i}{x-i}}.
\end{equation*}

\item \label{case:a=0b=0c<0d>0e=0} $a=0,b=0,c<0,d>0,e=0\Rightarrow \ga\ge 0, \gb\ge0, \gc=c=-l, \gd=0,\gg\ge0$. Then
\begin{equation*}
 (x+1)^l -C x^\ga (x-1)^\gb  (x+i)^\gg=D (x-i)^d.
\end{equation*}
Applying $x\Rightarrow -x$ and taking complex conjugation of Case \ref{case:a=0b>0c=0d<0e=0}) we get
\begin{equation*}
\boxed{f=\frac{1+ix}{1+x}, \quad \frac{x-i}{x+1}, \quad \frac{(1+i)(x-i)}{x+1}, \quad \frac{(1+i)(x-i)}{2(x+1)}}.
\end{equation*}

\item \label{case:a=0b=0c>0d<0e=0} $a=0,b=0,c>0,d<0,e=0\Rightarrow \ga\ge 0, \gb\ge0, \gc=0, \gd=d=-m,\gg\ge0$. Then
\begin{equation*}
   (x-i)^m -C x^\ga (x-1)^\gb  (x+i)^\gg=D (x+1)^c.
\end{equation*}
By Case \ref{case:a=0b=0c<0d>0e=0}) we get $\boxed{f=\frac{x+1}{i(x-i)}, \quad  \frac{x+1}{x-i}, \quad \frac{(1-i)(x+1)}{2(x-i)}, \quad \frac{(1-i)(x+1)}{x-i}}.$

\bigskip
\noindent\ \hskip-1.1cm \textbf{VII}. $a<0,b,c,d\ge0,e=0$.

\item \label{case:a<0b=0c=0d=0e=0} $a<0,b=c=d=e=0\Rightarrow \ga=a=-h, \gb\ge0,\gc\ge0,\gg\ge0$. Then
\begin{equation*}
  x^h-C (x-1)^\gb (x+1)^\gc (x-i)^\gd (x+i)^\gg=D.
\end{equation*}
By Case \ref{case:a>0b=0c=0d=0e=0}) we have
\begin{equation*}
\boxed{f=\frac{1}{x},\quad\frac{-1}{x},\quad\frac{i}{x},\quad\frac{-i}{x},\quad\frac{-1}{x^2},\quad\frac{1}{x^2},\quad\frac{1}{x^4}.} \end{equation*}

\item  \label{case:a<0b>0c=0d=0e=0} $a<0,b>0,c=d=e=0\Rightarrow \ga=a=-h, \gb=0, \gc\ge 0, \gd\ge 0,\gg\ge0$. Then
\begin{equation*}
  x^h-C (x+1)^\gc(x-i)^\gd (x+i)^\gg=D (x-1)^b.
\end{equation*}
By Case \ref{case:a>0b<0c=0d=0e=0}) we have
\begin{equation*}
\boxed{f=\frac{x-1}{x},\quad\frac{x-1}{2x},\quad \frac{(x-1)^2}{-4x},\quad\frac{x-1}{(1+i)x},\quad\frac{x-1}{(1-i)x},\quad\frac{(x-1)^2}{-2x}.}
\end{equation*}

\item  \label{case:a<0b=0c>0d=0e=0} $a<0,b=0,c>0,d=e=0\Rightarrow \ga=a=-h, \gb\ge0, \gc=0, \gd\ge0,\gg\ge0$. Then
\begin{equation*}
  x^h-C (x-1)^\gb(x-i)^\gd (x+i)^\gg=D (x+1)^c.
\end{equation*}
By Case \ref{case:a>0b=0c<0d=0e=0}) we get
\begin{equation*}
\boxed{f=\frac{x+1}{x},\quad  \frac{x+1}{2x},\quad \frac{(x+1)^2}{4x},\quad \frac{x+1}{(1+i)x},\quad\frac{x+1}{(1-i)x},\quad \frac{(x+1)^2}{2x}.}
\end{equation*}

\item \label{case:a<0b=0c=0d>0e=0} $a<0,b=c=0, d>0,e=0\Rightarrow \ga=a=-h, \gb\ge0, \gc\ge0, \gd=0,\gg\ge0$. Then
\begin{equation*}
  x^h-C (x-1)^\gb (x+1)^\gc (x+i)^\gg=D (x-i)^d.
\end{equation*}
By Case \ref{case:a>0b=0c=0d<0e=0}) we get
\begin{equation*}
\boxed{f=\frac{x-i}{x},\quad  \frac{x-i}{2x},\quad \frac{(x-i)^2}{-4ix},\quad \frac{x-i}{(1+i)x},\quad\frac{x-i}{(1-i)x},\quad \frac{(x-i)^2}{-2ix}.}
\end{equation*}

\item \label{case:a<0b>0c>0d=0e=0} $a<0,b>0,c>0,d=e=0 \Rightarrow \ga=a=-h,\gb=\gc=0, \gd\ge 0,\gg\ge0$.
\begin{equation*}
  x^h-C (x-i)^\gd (x+i)^\gg=D (x-1)^b (x+1)^c.
\end{equation*}
By Case \ref{case:a>0b<0c<0d=0e=0}) we get
\begin{equation*}
\boxed{f=\frac{x^2-1}{x^2}, \quad \frac{x^2-1}{2ix}, \quad \frac{x^2-1}{-2ix}, \quad \frac{x^2-1}{2x^2}, \quad \frac{(x^2-1)^2}{-4x^2}}.
\end{equation*}

\item \label{case:a<0b>0c=0d>0e=0} $a<0,b>0,c=0,d>0,e=0\Rightarrow \ga=a=-h,\gb=\gd=0, \gc,\gg\ge0$.
\begin{equation*}
  x^h -C (x+1)^\gc (x+i)^gg=D (x-1)^b (x-i)^d .
\end{equation*}
No solution exists by Case \ref{case:a>0b<0c=0d<0e=0} we get $\boxed{f=\frac{(x-1)(x-i)}{-2(1+i)x}}.$

\item \label{case:a<0b=0c>0d>0e=0} $a<0,b=0,c>0,d>0,e=0 \Rightarrow \ga=a=-h, \gb\ge 0, \gc=\gd=0, \gg\ge0$.
\begin{equation*}
  x^h-C (x-1)^\gb  (x+i)^\gg=D (x+1)^c(x-i)^d.
\end{equation*}
By Case \ref{case:a>0b=0c<0d<0e=0}) we get $\boxed{f=\frac{(x+1)(x-i)}{2(1-i)x}}.$

\item \label{case:a<0b>0c>0d>0e=0} $a<0,b>0,c>0,d>0,e=0 \Rightarrow \ga=a=-h,\gb=\gc=\gd=0,\gg\ge0$.
\begin{equation*}
  x^h-C  (x+i)^\gg=D (x-1)^b (x+1)^c(x-i)^d.
\end{equation*}
No solution exists by Case \ref{case:a>0b<0c<0d<0e=0}).

\bigskip
\noindent\ \hskip-1.1cm \textbf{VIII}. $a<0,b<0,e=0$.

\item  \label{case:a<0b<0c=0d=0e=0} $a<0,b<0,c=d=e=0\Rightarrow \ga=a=-h, \gb=b=-k, \gc\ge0, \gd\ge0,\gg\ge0$. Then
\begin{equation*}
 x^h (x-1)^k-C (x+1)^\gc(x-i)^\gd  (x+i)^\gg=D .
\end{equation*}
This is impossible by the 2-adic argument if we consider $x=0$ and $x=1$.

\item \label{case:a<0b<0c=0d>0e=0} $a<0,b<0,c=0,d>0,e=0\Rightarrow \ga=a=-h, \gb=b=-k, \gc\ge0,\gd=0,\gg\ge0$. Then
\begin{equation*}
 x^h (x-1)^k-C (x+1)^\gc  (x+i)^\gg=D (x-i)^d.
\end{equation*}
By Case \ref{case:a>0b>0c=0d<0e=0}) there is no solution.

\item \label{case:a<0b<0c>0d=0e=0}  $a<0,b<0,c>0,d=e=0\Rightarrow \ga=a=-h, \gb=b=-k, \gc=0,\gd\ge0,\gg\ge0$. Then
\begin{equation*}
 x^h (x-1)^k-C (x-i)^\gd (x+i)^\gg=D (x+1)^c .
\end{equation*}
By Case \ref{case:a>0b>0c<0d=0e=0}) we get $\boxed{f=\frac{x+1}{x(1-x)}}.$

\item \label{case:a<0b<0c>0d>0e=0} $a<0,b<0,c>0,d>0,e=0\Rightarrow \ga=a=-h, \gb=b=-k, \gc=\gd=0,\gg\ge0$. Then
\begin{equation*}
 x^h (x-1)^k-C (x+i)^\gg=D (x+1)^c (x-i)^d.
\end{equation*}
By Case \ref{case:a>0b>0c<0d<0e=0}) there is no solution.

\item \label{case:a<0b<0c<0d=0e=0} $a<0,b<0,c<0,d=e=0\Rightarrow \ga=a=-h, \gb=b=-k, \gc=c=-l, \gd\ge0,\gg\ge0$. Then
\begin{equation*}
 x^h (x-1)^k (x+1)^l-C(x-i)^\gd(x+i)^\gg=D  .
\end{equation*}
By Case \ref{case:a>0b>0c>0d=0e=0}) there is no solution.

\item \label{case:a<0b<0c=0d<0e=0} $a<0,b<0,c=0,d<0,e=0\Rightarrow \ga=a=-h, \gb=b=-k, \gc\ge0,\gd=d=-m, \gg\ge0$. Then
\begin{equation*}
 x^h (x-1)^k (x+1)^l-C (x-i)^\gd (x+i)^\gg=D.
\end{equation*}
By Case \ref{case:a>0b>0c=0d>0e=0}) there is no solution.

\item \label{case:a<0b<0c<0d>0e=0} $a<0,b<0,c<0,d>0,e=0\Rightarrow \ga=a=-h, \gb=b=-k, \gc=c=-l, \gd=0,\gg\ge0$. Then
\begin{equation*}
 x^h (x-1)^k (x+1)^l-C  (x+i)^\gg=D (x-i)^d.
\end{equation*}
By Case \ref{case:a>0b>0c>0d<0e=0}) we get $\boxed{f=\frac{i(x-i)^4}{8x(x^2-1)}}.$

\item \label{case:a<0b<0c>0d<0e=0} $a<0,b<0,c>0,d<0,e=0\Rightarrow \ga=a=-h, \gb=b=-k, \gc=0, \gd=d=-m,\gg\ge0$. Then
\begin{equation*}
 x^h (x-1)^k (x-i)^m-C  (x+i)^\gg=D  (x+1)^c.
\end{equation*}
By Case \ref{case:a>0b>0c<0d>0e=0}) there is no solution.

\item \label{case:a<0b<0c<0d<0e=0}  $a<0,b<0,c<0,d<0,e=0\Rightarrow \ga=a=-h, \gb=b=-k, \gc=c=-l, \gd=d=-m,\gg\ge0$. Then
\begin{equation*}
 x^h (x-1)^k (x+1)^l (x-i)^m-C  (x+i)^\gg=D
\end{equation*}
By Case \ref{case:a>0b>0c>0d>0e=0}) there is no solution.

\bigskip
\noindent\ \hskip-1.1cm \textbf{IX}. $a<0,b\ge 0,e=0,c<0$ or $d<0$.

\item \label{case:a<0b>0c<0d<0e=0} $a<0,b>0,c<0,d<0,e=0\Rightarrow \ga=a=-h, \gb=0, \gc=c=-l, \gd=d=-m,\gg\ge0$. Then
\begin{equation*}
 x^h (x+1)^l(x-i)^m-C  (x+i)^\gg=D (x-1)^b.
\end{equation*}
By Case \ref{case:a>0b<0c>0d>0e=0}) there is no solution.

\item \label{case:a<0b>0c<0d=0e=0} $a<0,b>0,c<0,d=e=0\Rightarrow \ga=a=-h, \gb=0, \gc=c=-l, \gd\ge 0,\gg\ge0$. Then
\begin{equation*}
 x^h (x+1)^l-C(x-i)^\gd  (x+i)^\gg=D (x-1)^b.
\end{equation*}
By Case \ref{case:a>0b<0c>0d=0e=0}) we get $\boxed{f=\frac{x-1}{x(x+1)}}.$

\item \label{case:a<0b>0c=0d<0e=0} $a<0,b>0,c=0,d<0,e=0\Rightarrow \ga=a=-h, \gb=0, \gc\ge 0, \gd=d=-m,\gg\ge0$. Then
\begin{equation*}
 x^h (x-i)^m -C (x+1)^\gc (x+i)^\gg=D (x-1)^b.
\end{equation*}
By Case \ref{case:a>0b<0c=0d>0e=0}) there is no solution.

\item \label{case:a<0b>0c<0d>0e=0} $a<0,b>0,c<0,d>0,e=0\Rightarrow \ga=a=-h, \gb=0, \gc=c=-l, \gd=0,\gg\ge0$. Then
\begin{equation*}
 x^h (x+1)^l-C  (x+i)^\gg=D (x-1)^b (x-i)^d.
\end{equation*}
By Case \ref{case:a>0b<0c>0d<0e=0}) there is no solution.

\item \label{case:a<0b>0c>0d<0e=0} $a<0,b>0,c>0,d<0,e=0\Rightarrow \ga=a=-h, \gb=0, \gc=0, \gd=d=-m,\gg\ge0$. Then
\begin{equation*}
 x^h (x-i)^m -C  (x+i)^\gg=D (x-1)^b (x+1)^c.
\end{equation*}
By Case \ref{case:a>0b<0c<0d>0e=0}) we get $\boxed{f=\frac{x^2-1}{x(x-i)}}$

\item \label{case:a<0b=0c<0d<0e=0} $a<0,b=0,c<0,d<0,e=0\Rightarrow \ga=a=-h, \gb\ge0, \gc=c=-l, \gd=d=-m,\gg\ge0$. Then
\begin{equation*}
 x^h (x+1)^l (x-i)^m-C (x-1)^\gb (x+i)^\gg=D.
\end{equation*}
By Case \ref{case:a>0b=0c>0d>0e=0}) there is no solution.

\item \label{case:a<0b=0c<0d=0e=0} $a<0,b=0,c<0,d=e=0\Rightarrow \ga=a=-h, \gb\ge0, \gc=c=-l, \gd\ge 0,\gg\ge0$. Then
\begin{equation*}
 x^h (x+1)^l-C (x-1)^\gb (x-i)^\gd (x+i)^\gg=D.
\end{equation*}
By Case \ref{case:a>0b=0c>0d=0e=0}) there is no solution.

\item \label{case:a<0b=0c=0d<0e=0} $a<0,b=0,c=0,d<0,e=0\Rightarrow \ga=a=-h, \gb\ge0, \gc\ge0, \gd=d=-m,\gg\ge0$. Then
\begin{equation*}
 x^h (x-i)^m-C (x-1)^\gb (x+1)^\gc (x+i)^\gg=D.
\end{equation*}
By Case \ref{case:a>0b=0c=0d>0e=0}) there is no solution.

\item \label{case:a<0b=0c<0d>0e=0} $a<0,b=0,c<0,d>0,e=0\Rightarrow \ga=a=-h, \gb\ge0, \gc=c=-l, \gd=0,\gg\ge0$. Then
\begin{equation*}
 x^h (x+1)^l -C (x-1)^\gb  (x+i)^\gg=D (x-i)^d.
\end{equation*}
By Case \ref{case:a>0b=0c>0d<0e=0}) there is no solution.

\item  \label{case:a<0b=0c>0d<0e=0} $a<0,b=0,c>0,d<0,e=0\Rightarrow \ga=a=-h, \gb\ge0, \gc=0, \gd=d=-m,\gg\ge0$. Then
\begin{equation*}
 x^h (x-i)^m-C (x-1)^\gb  (x+i)^\gg=D(x+1)^c.
\end{equation*}
By Case \ref{case:a>0b=0c<0d>0e=0}) there is no solution.

\bigskip
\noindent\ \hskip-1.1cm \textbf{X}. $a>0,b,c,d\ge0,e>0$.
There are eight cases none of which has a solution by the 2-adic argument if we consider $x=0$ and $x=-i$.

\newcommand\setItemnumber[1]{\setcounter{enumi}{\numexpr#1-1\relax}}

\setItemnumber{90}

In the following, we use the key idea that
\begin{itemize}
  \item Under $x\rightarrow -x$, $x-1 \leftrightharpoons x+1, x+i \leftrightharpoons x-1 $
  where $\leftharpoonup$ and $\rightharpoondown$ means the maps hold up to a constant. By our convention, this means $b\leftrightharpoons c$, $d \leftrightharpoons e$.

  \item Under $x\rightarrow -ix$, $x-1 \rightharpoondown x-i \rightharpoondown x+1 \rightharpoondown x+i \rightharpoondown x-1$,
  which means $b \rightharpoondown d \rightharpoondown c\rightharpoondown e \rightharpoondown b$.

  \item Under $x\rightarrow ix$, $x-1 \rightharpoondown x+i \rightharpoondown x+1 \rightharpoondown x-i \rightharpoondown x-1$,
  which means $b \rightharpoondown e \rightharpoondown c\rightharpoondown d \rightharpoondown b$.
\end{itemize}

\bigskip
\noindent\ \hskip-1.1cm \textbf{XI}. $a>0,b<0,e>0$.

\item \label{case:a>0b<0c=0d=0e>0} $a>0,b<0,c=0,d=0,e>0\Rightarrow \ga=0, \gb=b=-k, \gc\ge0, \gd\ge0, \gg\ge0$. Then
\begin{equation*}
 (x-1)^k-C (x+1)^\gc(x-i)^\gd=D x^a(x+i)^e.
\end{equation*}
Taking complex conjugation of all the coefficients in Case \ref{case:a>0b<0c=0d>0e=0}) we see there is no solution.

\item \label{case:a>0b<0c=0d>0e>0} $a>0,b<0,c=0,d>0,e>0\Rightarrow \ga=0, \gb=b=-k, \gc\ge0,\gd=0,\gg\ge0$. Then
\begin{equation*}
 (x-1)^k-C (x+1)^\gc=D x^a (x-i)^d(x+i)^e.
\end{equation*}
Applying $x\to ix$ in Case \ref{case:a>0b>0c>0d<0e=0}) we get $\boxed{f=\frac{-8x(x^2+1)}{(x-1)^4}}.$

\item\label{case:a>0b<0c>0d=0e>0} $a>0,b<0,c>0,d=0,e>0\Rightarrow \ga=\gc=0, \gb=b=-k, \gd\ge0, \gg\ge0$. Then
\begin{equation*}
 (x-1)^k-C (x-i)^\gd=D x^a (x+1)^c(x+i)^e.
\end{equation*}
Taking complex conjugation of all the coefficients in Case \ref{case:a>0b<0c>0d>0e=0}) we see that no solution exists.

\item\label{case:a>0b<0c>0d>0e>0} $a>0,b<0,c>0,d>0,e>0\Rightarrow \ga=0, \gb=b=-k, \gc=0, \gd=0,\gg\ge0$. Then
\begin{equation*}
 (x-1)^k-C=D x^a (x+1)^c (x-i)^d(x+i)^e.
\end{equation*}
No solution exists by Case \ref{case:a=0b>0c=0d=0e=0}.

\item \label{case:a>0b<0c<0d=0e>0} $a>0,b<0,c<0,d=0,e>0\Rightarrow \ga=0, \gb=b=-k, \gc=c=-l, \gd\ge0,\gg\ge0$. Then
\begin{equation*}
 (x-1)^k (x+1)^l-C (x-i)^\gd=D x^a(x+i)^e.
\end{equation*}
Taking complex conjugation of all the coefficients in Case \ref{case:a>0b<0c<0d>0e=0}) we get $\boxed{f=\frac{x(x+i)}{x^2-1}}$

\item\label{case:a>0b<0c=0d<0e>0} $a>0,b<0,c=0,d<0,e>0\Rightarrow \ga=0, \gb=b=-k, \gc\ge0, \gd=d=-m,\gg\ge0$. Then
\begin{equation*}
 (x-1)^k(x-i)^m -C(x+1)^\gc=D x^a(x+i)^e.
\end{equation*}
Applying $x\to ix$ in Case \ref{case:a>0b>0c<0d<0e=0}) we see that no solution exists.

\item\label{case:a>0b<0c<0d>0e>0} $a>0,b<0,c<0,d>0,e>0\Rightarrow \ga=0, \gb=b=-k, \gc=c=-l, \gd=0,\gg\ge0$. Then
\begin{equation*}
 (x-1)^k (x+1)^l-C=D x^a(x-i)^d(x+i)^e.
\end{equation*}
No solution exists by Case \ref{case:a=0b>0c>0d=0e=0}.

\item \label{case:a>0b<0c>0d<0e>0} $a>0,b<0,c>0,d<0,e>0\Rightarrow \ga=\gc=0, \gb=b=-k, \gd=d=m, \gg\ge0$. Then
\begin{equation*}
   (x-1)^k (x-i)^m -C=D x^a (x+1)^c(x+i)^e.
\end{equation*}
No solution exists by Case \ref{case:a=0b>0c=0d>0e=0}.

\item\label{case:a>0b<0c<0d<0e>0} $a>0,b<0,c<0,d<0,e>0\Rightarrow \ga=0, \gb=b=-k, \gc=c=-l, \gd=d=-m,\gg\ge0$. Then
\begin{equation*}
 (x-1)^k (x+1)^l(x-i)^m-C=D x^a(x+i)^e.
\end{equation*}
No solution exists by Case \ref{case:a=0b>0c>0d>0e=0}.

\bigskip
\bigskip\noindent\ \hskip-1.1cm \textbf{XII}. $a>0,b\ge0,e>0,c<0$ or $d<0$.

\item \label{case:a>0b>0c<0d<0e>0} $a>0,b>0,c<0,d<0,e>0\Rightarrow \ga=\gb=0, \gc=c=-l, \gd=d=-m,\gg\ge0$. Then
\begin{equation*}
 (x+1)^l(x-i)^m-C=D x^a(x-1)^b(x+i)^e.
\end{equation*}
No solution exists by Case \ref{case:a=0b=0c>0d>0e=0}.

\item \label{case:a>0b>0c<0d=0e>0} $a>0,b>0,c<0,d=0,e>0\Rightarrow \ga=0, \gb=0, \gc=c=-l, \gd\ge 0,\gg\ge0$. Then
\begin{equation*}
 (x+1)^l-C(x-i)^\gd=D x^a(x-1)^b(x+i)^e.
\end{equation*}
Taking complex conjugation of all the coefficients in Case \ref{case:a>0b>0c<0d>0e=0}) we see that no solution exists.

\item \label{case:a>0b>0c=0d<0e>0} $a>0,b>0,c=0,d<0,e>0\Rightarrow \ga=\gb=0,\gc\ge 0,\gd=d=-m,\gg\ge0$. Then
\begin{equation*}
 (x-i)^m-C (x+1)^\gc=D x^a(x-1)^b(x+i)^e.
\end{equation*}
Applying $x\to ix$ in Case \ref{case:a>0b>0c<0d>0e=0}) we see that no solution exists.

\item \label{case:a>0b>0c<0d>0e>0} $a>0,b>0,c<0,d>0,e>0\Rightarrow \ga=0, \gb=0, \gc=c=-l, \gd=0,\gg\ge0$. Then
\begin{equation*}
 (x+1)^l-C=D x^a(x-1)^b (x-i)^d(x+i)^e.
\end{equation*}
By Case \ref{case:a=0b=0c>0d=0e=0}) we see that no solution exists.

\item \label{case:a>0b>0c>0d<0e>0} $a>0,b>0,c>0,d<0,e>0\Rightarrow \ga=\gb=\gc=0, \gd=d=-m,\gg\ge0$. Then
\begin{equation*}
(x-i)^m-C=D x^a(x-1)^b (x+1)^c (x+i)^e.
\end{equation*}
By Case \ref{case:a=0b=0c=0d>0e=0}) we see that no solution exists.

\item \label{case:a>0b=0c<0d<0e>0} $a>0,b=0,c<0,d<0,e>0\Rightarrow \ga=0, \gb\ge0, \gc=c=-l, \gd=d=-m,\gg\ge0$. Then
\begin{equation*}
 (x+1)^l (x-i)^m-C (x-1)^\gb=D x^a(x+i)^e.
\end{equation*}
Applying $x\to -ix$ in Case \ref{case:a>0b<0c>0d<0e=0}) we see that no solution exists.

\item \label{case:a>0b=0c<0d=0e>0} $a>0,b=0,c<0,d=0,e>0\Rightarrow \ga=0, \gb\ge0, \gc=c=-l, \gd\ge 0,\gg\ge0$. Then
\begin{equation*}
 (x+1)^l-C (x-1)^\gb (x-i)^\gd=D x^a(x+i)^e.
\end{equation*}
Taking complex conjugation of all the coefficients in Case \ref{case:a>0b=0c<0d>0e=0}) we see that no solution exists.

\item \label{case:a>0b=0c=0d<0e>0} $a>0,b=0,c=0,d<0,e>0\Rightarrow \ga=0,\gb\ge0,\gc\ge0, \gd=d=-m,\gg\ge0$. Then
\begin{equation*}
 (x-i)^m-C (x-1)^\gb(x+1)^\gc=D x^a (x+i)^e.
\end{equation*}
Applying $x\to -ix$ in Case \ref{case:a>0b<0c>0d=0e=0}) we get $\boxed{f=\frac{-ix(x+i)}{x-i}}.$

\item \label{case:a>0b=0c<0d>0e>0} $a>0,b=0,c<0,d>0,e>0\Rightarrow \ga=0, \gb\ge0, \gc=c=-l, \gd=0,\gg\ge0$. Then
\begin{equation*}
 (x+1)^l -C (x-1)^\gb=D x^a (x-i)^d(x+i)^e.
\end{equation*}
Applying $x\to -ix$ in Case \ref{case:a>0b>0c>0d<0e=0}) we get $\boxed{f=\frac{8x(x^2+1)}{(x+1)^4}}.$

\item \label{case:a>0b=0c>0d<0e>0} $a>0,b=0,c>0,d<0,e>0\Rightarrow \ga=0, \gb\ge0, \gc=0, \gd=d=-m,\gg\ge0$. Then
\begin{equation*}
 (x-i)^m -C (x-1)^\gb=D x^a (x+1)^c(x+i)^e.
\end{equation*}
Applying $x\to -ix$ in Case \ref{case:a>0b>0c<0d>0e=0}) we see that no solution exists.

\newpage
\noindent\ \hskip-1.1cm \textbf{XIII}. $a=0,b,c,d\ge0,e>0$.

\item \label{case:a=0b=0c=0d=0e>0} $a=0,b=c=d=0,e>0\Rightarrow \ga,\gb,\gc,\gd\ge0,\gg=0$.  Then
\begin{equation*}
 1-C x^\ga (x-1)^\gb (x+1)^\gc (x-i)^gd=D (x+i)^e.
\end{equation*}
Taking complex conjugation of all the coefficients in Case \ref{case:a=0b=0c=0d>0e=0}) we get
\begin{equation*}
\boxed{f=\frac{x+i}{i}, \quad\frac{x+i}{2i}, \quad \frac{x+i}{i-1}, \quad \frac{x+i}{1+i}}.
\end{equation*}

\item \label{case:a=0b>0c=0d=0e>0} $a=0,b>0,c=d=0,e>0\Rightarrow \ga\ge 0,\gb=0, \gc\ge 0, \gd\ge 0,\gg\ge0$. Then we get
\begin{equation*}
 B-x^\ga (x+1)^\gc (x-i)^\gd=-(x-1)^b(x+i)^e.
\end{equation*}
Taking complex conjugation of all the coefficients in Case \ref{case:a=0b>0c=0d>0e=0}) we see that no solution exists.

\item \label{case:a=0b=0c>0d=0e>0} $a=b=0,c>0,d=0,e>0\Rightarrow \ga\ge0, \gb\ge0,\gc=0, \gd\ge0,\gg\ge0$. Then we must have
\begin{equation*}
 1-Cx^\ga (x-1)^\gb(x-i)^\gd=D(x+1)^c(x+i)^e.
\end{equation*}
Taking complex conjugation of all the coefficients in Case \ref{case:a=0b=0c>0d>0e=0}) we see that no solution exists.

\item \label{case:a=0b=0c=0d>0e>0} $a=b=c=0,d>0,e>0\Rightarrow \ga\ge0, \gb\ge0,\gc\ge0, \gd=0,\gg\ge0$. By symmetry $i\leftrightarrow -i$ from above we get
\begin{equation*}
  1-Cx^\ga (x-1)^\gb (x+1)^\gc=D(x-i)^d(x+i)^e.
\end{equation*}
Applying $x\to ix$ in Case \ref{case:a=0b>0c>0d=0e=0}) we get $\boxed{f=x^2+1}$ or $\boxed{f=\frac{1+x^2}2}.$

\item \label{case:a=0b>0c>0d=0e>0} $a=0,b>0,c>0, d=0,e>0\Rightarrow \ga=\gb=\gc=0, \gd\ge 0,\gg\ge0$. Then we must have
\begin{equation*}
  B-x^\ga (x-i)^\gd=-(x-1)^b (x+1)^c(x+i)^e.
\end{equation*}
Taking complex conjugation of all the coefficients in Case \ref{case:a=0b>0c>0d>0e=0}) we see that no solution exists.

\item \label{case:a=0b>0c=0d>0e>0} $a=0,b>0,c=0,d>0,e>0\Rightarrow \ga\ge0, \gb=0, \gc\ge0, \gd=0,\gg\ge0$. Then we must have
\begin{equation*}
  B-x^\ga (x+1)^\gc=-(x-1)^b (x-i)^d(x+i)^e.
\end{equation*}
Applying $x\to ix$ in Case \ref{case:a=0b>0c>0d>0e=0}) we see that no solution exists.

\item \label{case:a=0b=0c>0d>0e>0} $a=0,b=0,c>0, d>0,e>0\Rightarrow \ga\ge 0, \gb\ge 0, \gc=\gd=0, \gg\ge0$. Then we must have
\begin{equation*}
  1-Cx^\ga (x-1)^\gb=D(x+1)^c(x-i)^d(x+i)^e.
\end{equation*}
Applying $x\to ix$ in Case \ref{case:a=0b>0c>0d>0e=0}) we see that no solution exists.

\item \label{case:a=0b>0c>0d>0e>0} $a=0,b>0,c>0,d>0,e>0\Rightarrow \ga=\gc=\gb=0, \gg\ge0$. Then we must have
\begin{equation*}
  B- x^\ga=-(x-1)^b (x+1)^c (x-i)^d(x+i)^e.
\end{equation*}
By Case \ref{case:a>0b=0c=0d=0e=0}) we get $\boxed{f=1-x^4}.$

\newpage\noindent\ \hskip-1.1cm \textbf{XIV}. $a=0,b<0,e>0$.

\item \label{case:a=0b<0c=0d=0e>0} $a=0,b<0,c=0,d=0,e>0\Rightarrow \ga\ge 0, \gb=b=-k, \gc\ge0, \gd\ge0,\gg\ge0$. Then
\begin{equation*}
 (x-1)^k-C x^\ga (x+1)^\gc (x-i)^\gd=D (x+i)^e.
\end{equation*}
Taking complex conjugation of all the coefficients in Case \ref{case:a=0b<0c=0d>0e=0}) we get
\begin{equation*}
\boxed{f=\frac{1-ix}{1-x}, \quad  \frac{x+i}{x-1}, \quad \frac{(1+i)(x+i)}{x-1}, \quad\frac{(1+i)(x+i)}{2(x-1)}}.
\end{equation*}

\item \label{case:a=0b<0c=0d>0e>0} $a=0,b<0,c=0,d>0,e>0\Rightarrow \ga\ge 0, \gb=b=-k, \gc\ge0,\gd=0,\gg\ge0$. Then
\begin{equation*}
 (x-1)^k-C x^\ga (x+1)^\gc=D (x-i)^d(x+i)^e.
\end{equation*}
Applying $x\to ix$ in Case \ref{case:a=0b>0c>0d<0e=0}) we get
$\boxed{f=\frac{x^2+1}{(x-1)^2},\ \frac{2(x^2+1)}{(x-1)^2}, \ \frac{1+x^2}{1-x}}.$

\item \label{case:a=0b<0c>0d=0e>0} $a=0,b<0,c>0,d=0,e>0\Rightarrow \ga\ge 0, \gb=b=-k, \gc=0,\gd\ge0,\gg\ge0$. Then
\begin{equation*}
 (x-1)^k-C x^\ga(x-i)^\gd=D (x+1)^c(x+i)^e.
\end{equation*}
Taking complex conjugation of all the coefficients in Case \ref{case:a=0b<0c>0d>0e=0}) we see that no solution exists.

\item \label{case:a=0b<0c>0d>0e>0} $a=0,b<0,c>0,d>0,e>0\Rightarrow \ga\ge 0, \gb=b=-k, \gc=0, \gd=0,\gg\ge0$. Then
\begin{equation*}
 (x-1)^k-Cx^\ga=D (x+1)^c (x-i)^d(x+i)^e.
\end{equation*}

i) If $\ga=0$ then by Case \ref{case:a=0b>0c=0d=0e=0}) we see that no solution exists to make $c,d,e>0$.

ii) If $\ga>0$ then by Case \ref{case:a<0b>0c=0d=0e=0}) we see that no solution exists to make $c,d,e>0$.

\item \label{case:a=0b<0c<0d=0e>0} $a=0,b<0,c<0,d=0,e>0\Rightarrow \ga\ge 0, \gb=b=-k, \gc=c=-l, \gd\ge0,\gg\ge0$. Then
\begin{equation*}
 (x-1)^k (x+1)^l-C x^\ga (x-i)^\gd=D(x+i)^e.
\end{equation*}
Taking complex conjugation of all the coefficients in Case \ref{case:a=0b<0c<0d>0e=0}) we get
\begin{equation*}
\boxed{f=\frac{(x+i)^2}{x^2-1},\quad \frac{(x+i)^2}{2(x^2-1)},\quad \frac{ix-1}{x^2-1}}.
\end{equation*}

\item\label{case:a=0b<0c=0d<0e>0} $a=0,b<0,c=0,d<0,e>0\Rightarrow \ga\ge 0, \gb=b=-k, \gc\ge0, \gd=d=-m,\gg\ge0$. Then
\begin{equation*}
 (x-1)^k(x-i)^m -C x^\ga (x+1)^\gc=D(x+i)^e.
\end{equation*}
Applying $x\to ix$ in Case \ref{case:a=0b>0c<0d<0e=0}) we see that no solution exists.

\item \label{case:a=0b<0c<0d>0e>0} $a=0,b<0,c<0,d>0,e>0\Rightarrow \ga\ge 0, \gb=b=-k, \gc=c=-l, \gd=0,\gg\ge0$. Then
\begin{equation*}
 (x-1)^k (x+1)^l-C x^\ga=D (x-i)^d(x+i)^e.
\end{equation*}

i) $\ga=0$. By Case \ref{case:a=0b>0c>0d=0e=0}) we get $\boxed{f=\frac{x^2+1}{x^2-1}, \frac{x^2+1}{1-x^2}}.$

ii) $\ga\ge 0$.  By Case \ref{case:a<0b>0c>0d=0e=0}) we get $\boxed{f=\frac{(x^2+1)^2}{(x^2-1)^2}}.$

\item \label{case:a=0b<0c>0d<0e>0} $a=0,b<0,c>0,d<0,e>0\Rightarrow \ga\ge 0, \gb=b=-k, \gc=0, \gd=d=-m,\gg\ge0$. Then
\begin{equation*}
 (x-1)^k (x-i)^m -C x^\ga=D (x+1)^c(x+i)^e.
\end{equation*}

i) $\ga=0$. By Case \ref{case:a=0b>0c=0d>0e=0}) we see that no solution exists.

ii) $\ga\ge 0$. By Case \ref{case:a<0b>0c=0d>0e=0}) we get $\boxed{f=\frac{(x+1)(x+i)}{(x-1)(x-i)}}.$

\item \label{case:a=0b<0c<0d<0e>0} $a=0,b<0,c<0,d<0,e>0\Rightarrow \ga\ge 0, \gb=b=-k, \gc=c=-l, \gd=d=-m,\gg\ge0$. Then
\begin{equation*}
 (x-1)^k (x+1)^l(x-i)^m-C x^\ga=D(x+i)^e.
\end{equation*}

i) $\ga=0$. By Case \ref{case:a=0b>0c>0d>0e=0}) we see that no solution exists.

ii) $\ga\ge 0$. By Case \ref{case:a<0b>0c>0d>0e=0}) we see that no solution exists.

\bigskip\noindent\ \hskip-1.1cm \textbf{XV}. $a=0,b\ge 0,e>0,c<0$ or $d<0$.

\item \label{case:a=0b>0c<0d<0e>0} $a=0,b>0,c<0,d<0,e>0\Rightarrow \ga\ge 0, \gb=0, \gc=c=-l, \gd=d=-m,\gg\ge0$. Then
\begin{equation*}
 (x+1)^l(x-i)^m-C x^\ga=D (x-1)^b(x+i)^e.
\end{equation*}

i) $\ga=0$. By Case \ref{case:a=0b=0c>0d>0e=0}) we see that no solution exists.

ii) $\ga\ge 0$. By Case \ref{case:a<0b=0c>0d>0e=0}) we get $\boxed{f=\frac{(x-1)(x+i)}{(x+1)(x-i)}}.$

\item \label{case:a=0b>0c<0d=0e>0} $a=0,b>0,c<0,d=0,e>0\Rightarrow \ga\ge 0, \gb=0, \gc=c=-l, \gd\ge 0,\gg\ge0$. Then
\begin{equation*}
 (x+1)^l-Cx^\ga (x-i)^\gd=D (x-1)^b(x+i)^e.
\end{equation*}
Taking complex conjugation of all the coefficients in Case \ref{case:a=0b>0c<0d>0e=0}) we see that no solution exists.

\item \label{case:a=0b>0c=0d<0e>0} $a=0,b>0,c=0,d<0,e>0\Rightarrow \ga\ge 0, \gb=0, \gc\ge 0, \gd=d=-m,\gg\ge0$. Then
\begin{equation*}
 (x-i)^m-Cx^\ga (x+1)^\gc=D (x-1)^b(x+i)^e.
\end{equation*}
Applying $x\to ix$ in Case \ref{case:a=0b>0c<0d>0e=0}) we see that no solution exists.

\item \label{case:a=0b>0c<0d>0e>0} $a=0,b>0,c<0,d>0,e>0\Rightarrow \ga\ge 0, \gb=0, \gc=c=-l, \gd=0,\gg\ge0$. Then
\begin{equation*}
 (x+1)^l-C x^\ga=D (x-1)^b (x-i)^d(x+i)^e.
\end{equation*}

i) $\ga=0$. By Case \ref{case:a=0b=0c>0d=0e=0}) we see that no solution exists.

ii) $\ga\ge 0$. By Case \ref{case:a<0b=0c>0d=0e=0}) we see that no solution exists.

\item \label{case:a=0b>0c>0d<0e>0} $a=0,b>0,c>0,d<0,e>0\Rightarrow \ga\ge 0, \gb=0, \gc=0, \gd=d=-m,\gg\ge0$. Then
\begin{equation*}
 (x-i)^m-C x^\ga=D (x-1)^b (x+1)^c(x+i)^e.
\end{equation*}

i) $\ga=0$. By Case \ref{case:a=0b=0c=0d>0e=0}) we see that no solution exists.

ii) $\ga\ge 0$. By Case \ref{case:a<0b=0c=0d>0e=0}) we see that no solution exists.

\item \label{case:a=0b=0c<0d<0e>0} $a=0,b=0,c<0,d<0,e>0\Rightarrow \ga\ge 0, \gb\ge0, \gc=c=-l, \gd=d=-m,\gg\ge0$. Then
\begin{equation*}
 (x+1)^l (x-i)^m-C x^\ga (x-1)^\gb=D(x+i)^e.
\end{equation*}
Applying $x\to -ix$ in Case \ref{case:a=0b<0c>0d<0e=0}) we see that no solution exists.

\item \label{case:a=0b=0c<0d=0e>0} $a=0,b=0,c<0,d=0,e>0\Rightarrow \ga\ge 0, \gb\ge0, \gc=c=-l, \gd\ge 0,\gg\ge0$. Then
\begin{equation*}
 (x+1)^l-C x^\ga (x-1)^\gb (x-i)^\gd=D(x+i)^e.
\end{equation*}
Taking complex conjugation of all the coefficients in Case \ref{case:a=0b=0c<0d>0e=0}) we get
\begin{equation*}
\boxed{f=\frac{1-ix}{1+x}, \quad  \frac{x+i}{x+1}, \quad \frac{(1-i)(x+i)}{x+1}, \quad \frac{(1-i)(x+i)}{2(x+1)}}.
\end{equation*}

\item \label{case:a=0b=0c=0d<0e>0} $a=0,b=0,c=0,d<0,e>0\Rightarrow \ga\ge 0, \gb\ge0, \gc\ge 0, \gd=d=-m,\gg\ge0$. Then
\begin{equation*}
 (x-i)^m-C x^\ga (x+1)^\gc (x-1)^\gb=D(x+i)^e.
\end{equation*}
Applying $x\to ix$ in Case \ref{case:a=0b>0c<0d=0e=0}) we get
\begin{equation*}
\boxed{f=\frac{x+i}{x-i}, \quad \frac{x+i}{i-x}, \quad \frac{(x+i)^2}{(x-i)^2}, \quad \frac{x+i}{i(x-i)}, \quad \frac{x+i}{i(i-x)},  \quad \frac{-(x+i)^2}{(x-i)^2}, \quad  \frac{(x+i)^4}{(x-i)^4}}.
\end{equation*}

\item \label{case:a=0b=0c<0d>0e>0} $a=0,b=0,c<0,d>0,e>0\Rightarrow \ga\ge 0, \gb\ge0, \gc=c=-l, \gd=0,\gg\ge0$. Then
\begin{equation*}
 (x+1)^l -C x^\ga (x-1)^\gb=D (x-i)^d(x+i)^e.
\end{equation*}
Applying $x\to -ix$ in Case \ref{case:a=0b>0c>0d<0e=0}) we get
$\boxed{f=\frac{x^2+1}{(x+1)^2},\ \frac{2(x^2+1)}{(x+1)^2}, \ \frac{x^2+1}{x+1}}.$

\item \label{case:a=0b=0c>0d<0e>0} $a=0,b=0,c>0,d<0,e>0\Rightarrow \ga\ge 0, \gb\ge0, \gc=0, \gd=d=-m,\gg\ge0$. Then
\begin{equation*}
   (x-i)^m -C x^\ga (x-1)^\gb=D (x+1)^c(x+i)^e.
\end{equation*}
Applying $x\to -ix$ in Case \ref{case:a=0b<0c>0d>0e=0}) we see that no solution exists.

\bigskip\noindent\ \hskip-1.1cm \textbf{XVI}. $a<0,b,c,d\ge0,e>0$.

\item \label{case:a<0b=0c=0d=0e>0} $a<0,b=c=d=0,e>0\Rightarrow \ga=a=-h, \gb\ge0,\gc\ge0,\gg\ge0$. Then
\begin{equation*}
  x^h-C (x-1)^\gb (x+1)^\gc (x-i)^\gd=D(x+i)^e.
\end{equation*}
Taking complex conjugation of all the coefficients in Case \ref{case:a<0b=0c=0d>0e=0}) we get
\begin{equation*}
\boxed{f=\frac{x+i}{x},\quad  \frac{x+i}{2x},\quad \frac{(x+i)^2}{4ix},\quad \frac{x+i}{(1-i)x},\quad\frac{x+i}{(1+i)x},\quad \frac{(x+i)^2}{2ix}.}
\end{equation*}

\item  \label{case:a<0b>0c=0d=0e>0} $a<0,b>0,c=d=0,e>0\Rightarrow \ga=a=-h, \gb=0, \gc\ge 0, \gd\ge 0,\gg\ge0$. Then
\begin{equation*}
  x^h-C (x+1)^\gc(x-i)^\gd=D (x-1)^b(x+i)^e.
\end{equation*}
Taking complex conjugation of all the coefficients in Case \ref{case:a<0b>0c=0d>0e=0}) we get $\boxed{f=\frac{(x-1)(x+i)}{-2(1-i)x}}.$

\item  \label{case:a<0b=0c>0d=0e>0} $a<0,b=0,c>0,d=0,e>0\Rightarrow \ga=a=-h, \gb\ge0, \gc=0, \gd\ge0,\gg\ge0$. Then
\begin{equation*}
  x^h-C (x-1)^\gb(x-i)^\gd=D (x+1)^c(x+i)^e.
\end{equation*}
Taking complex conjugation of all the coefficients in Case \ref{case:a<0b=0c>0d>0e=0}) we get $\boxed{f=\frac{(x+1)(x+i)}{2(1+i)x}}.$

\item \label{case:a<0b=0c=0d>0e>0} $a<0,b=c=0, d>0,e>0\Rightarrow \ga=a=-h, \gb\ge0, \gc\ge0, \gd=0,\gg\ge0$. Then
\begin{equation*}
  x^h-C (x-1)^\gb (x+1)^\gc=D (x-i)^d(x+i)^e.
\end{equation*}
Applying $x\to ix$ in Case \ref{case:a<0b>0c>0d=0e=0}) we get
\begin{equation*}
\boxed{f=\frac{x^2+1}{x^2}, \quad \frac{x^2+1}{2x}, \quad \frac{x^2+1}{-2x}, \quad \frac{x^2+1}{2x^2}, \quad \frac{(x^2+1)^2}{4x^2}}.
\end{equation*}

\item \label{case:a<0b>0c>0d=0e>0} $a<0,b>0,c>0,d=0,e>0\Rightarrow \ga=a=-h,\gb=\gc=0, \gd\ge 0,\gg\ge0$.
\begin{equation*}
  x^h-C (x-i)^\gd=D (x-1)^b (x+1)^c(x+i)^e.
\end{equation*}
Taking complex conjugation of all the coefficients in Case \ref{case:a<0b>0c>0d>0e=0}) we see that no solution exists.

\item \label{case:a<0b>0c=0d>0e>0} $a>0,b>0,c=0,d>0,e>0\Rightarrow \ga=\gb=\gd=0, \gc\ge 0, \gg\ge0$.
\begin{equation*}
  1-C x^\ga (x+1)^\gc =D (x-1)^b (x-i)^d (x+i)^e.
\end{equation*}
Applying $x\to ix$ in Case \ref{case:a<0b>0c>0d>0e=0}) we see that no solution exists.

\item \label{case:a<0b=0c>0d>0e>0} $a<0,b=0,c>0,d>0,e>0\Rightarrow \ga=a=-h, \gb\ge 0, \gc=\gd=0, \gg\ge0$.
\begin{equation*}
  x^h-C (x-1)^\gb=D (x+1)^c(x-i)^d(x+i)^e.
\end{equation*}
Applying $x\to -ix$ in Case \ref{case:a<0b>0c>0d>0e=0}) we see that no solution exists.

\item \label{case:a<0b>0c>0d>0e>0} $a<0,b>0,c>0,d>0,e>0\Rightarrow \ga=a=-h,\gb=\gc=\gd=0,\gg\ge0$.
\begin{equation*}
  x^h-C=D (x-1)^b (x+1)^c(x-i)^d(x+i)^e.
\end{equation*}
By Case \ref{case:a>0b=0c=0d=0e=0}) we get $\boxed{f=\frac{x^4-1}{x^4}}.$

\bigskip\noindent\ \hskip-1.1cm \textbf{XVII}. $a<0,b<0,e>0$.

\item  \label{case:a<0b<0c=0d=0e>0} $a<0,b<0,c=d=0,e>0\Rightarrow \ga=a=-h, \gb=b=-k, \gc\ge0, \gd\ge0,\gg\ge0$. Then
\begin{equation*}
 x^h (x-1)^k-C (x+1)^\gc(x-i)^\gd=D (x+i)^e.
\end{equation*}
Taking complex conjugation of all the coefficients in Case \ref{case:a<0b<0c=0d>0e=0}) we see that no solution exists.

\item \label{case:a<0b<0c=0d>0e>0} $a<0,b<0,c=0,d>0,e>0\Rightarrow \ga=a=-h, \gb=b=-k, \gc\ge0,\gd=0,\gg\ge0$. Then
\begin{equation*}
 x^h (x-1)^k-C (x+1)^\gc=D (x-i)^d(x+i)^e.
\end{equation*}
Applying $x\to ix$ in Case \ref{case:a<0b>0c>0d<0e=0}) we get $\boxed{f=\frac{x^2+1}{x(x-1)}}$

\item \label{case:a<0b<0c>0d=0e>0}  $a<0,b<0,c>0,d=0,e>0\Rightarrow \ga=a=-h, \gb=b=-k, \gc=0,\gd\ge0,\gg\ge0$. Then
\begin{equation*}
 x^h (x-1)^k-C (x-i)^\gd=D (x+1)^c (x+i)^e.
\end{equation*}
Applying $x\to ix$ in Case \ref{case:a<0b<0c>0d>0e=0}) we see that no solution exists.

\item \label{case:a<0b<0c>0d>0e>0} $a<0,b<0,c>0,d>0,e>0\Rightarrow \ga=a=-h, \gb=b=-k, \gc=\gd=0,\gg\ge0$. Then
\begin{equation*}
 x^h (x-1)^k-C=D (x+1)^c (x-i)^d(x+i)^e.
\end{equation*}
By Case \ref{case:a>0b>0c=0d=0e=0}) we see that no solution exists.

\item \label{case:a<0b<0c<0d=0e>0} $a<0,b<0,c<0,d=0,e>0\Rightarrow \ga=a=-h, \gb=b=-k, \gc=c=-l, \gd\ge0,\gg\ge0$. Then
\begin{equation*}
 x^h (x-1)^k (x+1)^l-C(x-i)^\gd=D  (x+i)^e.
\end{equation*}
Taking complex conjugation of all the coefficients in Case \ref{case:a<0b<0c<0d>0e=0}) we get $\boxed{f=\frac{-i(x+i)^4}{8x(x^2-1)}}.$

\item \label{case:a<0b<0c=0d<0e>0} $a<0,b<0,c=0,d<0,e>0\Rightarrow \ga=a=-h, \gb=b=-k, \gc\ge0,\gd=d=-m, \gg\ge0$. Then
\begin{equation*}
 x^h (x-1)^k (x+1)^l-C (x-i)^\gd=D(x+i)^e.
\end{equation*}
Applying $x\to ix$ in Case \ref{case:a<0b>0c<0d<0e=0}) we see that no solution exists.

\item \label{case:a<0b<0c<0d>0e>0} $a<0,b<0,c<0,d>0,e>0\Rightarrow \ga=a=-h, \gb=b=-k, \gc=c=-l, \gd=0,\gg\ge0$. Then
\begin{equation*}
 x^h (x-1)^k (x+1)^l-C=D (x-i)^d(x+i)^e.
\end{equation*}
By Case \ref{case:a>0b>0c>0d=0e=0}) we see that no solution exists.

\item \label{case:a<0b<0c>0d<0e>0} $a<0,b<0,c>0,d<0,e>0\Rightarrow \ga=a=-h, \gb=b=-k, \gc=0, \gd=d=-m,\gg\ge0$. Then
\begin{equation*}
 x^h (x-1)^k (x-i)^m-C=D  (x+1)^c(x+i)^e.
\end{equation*}
By Case \ref{case:a>0b>0c=0d>0e=0}) we see that no solution exists.

\item \label{case:a<0b<0c<0d<0e>0}  $a<0,b<0,c<0,d<0,e>0\Rightarrow \ga=a=-h, \gb=b=-k, \gc=c=-l, \gd=d=-m,\gg\ge0$. Then
\begin{equation*}
 x^h (x-1)^k (x+1)^l (x-i)^m-C=D
\end{equation*}
By Case \ref{case:a>0b>0c>0d>0e=0}) we see that no solution exists.

\bigskip\noindent\ \hskip-1.1cm \textbf{XVIII}. $a<0,b\ge 0,e>0,c<0$ or $d<0$.

\item \label{case:a<0b>0c<0d<0e>0} $a<0,b>0,c<0,d<0,e>0\Rightarrow \ga=a=-h, \gb=0, \gc=c=-l, \gd=d=-m,\gg\ge0$. Then
\begin{equation*}
 x^h (x+1)^l(x-i)^m-C=D (x-1)^b(x+i)^e.
\end{equation*}
By Case \ref{case:a>0b=0c>0d>0e=0}) we see that no solution exists.

\item \label{case:a<0b>0c<0d=0e>0} $a<0,b>0,c<0,d=0,e>0\Rightarrow \ga=a=-h, \gb=0, \gc=c=-l, \gd\ge 0,\gg\ge0$. Then
\begin{equation*}
 x^h (x+1)^l-C(x-i)^\gd=D (x-1)^b(x+i)^e.
\end{equation*}
Taking complex conjugation of all the coefficients in Case \ref{case:a<0b>0c<0d>0e=0})

\item \label{case:a<0b>0c=0d<0e>0} $a<0,b>0,c=0,d<0,e>0\Rightarrow \ga=a=-h, \gb=0, \gc\ge 0, \gd=d=-m,\gg\ge0$. Then
\begin{equation*}
 x^h (x-i)^m -C (x+1)^\gc=D (x-1)^b(x+i)^e.
\end{equation*}
Applying $x\to ix$ in Case \ref{case:a<0b>0c<0d>0e=0}) we see that no solution exists.

\item \label{case:a<0b>0c<0d>0e>0} $a<0,b>0,c<0,d>0,e>0\Rightarrow \ga=a=-h, \gb=0, \gc=c=-l, \gd=0,\gg\ge0$. Then
\begin{equation*}
 x^h (x+1)^l-C=D (x-1)^b (x-i)^d(x+i)^e.
\end{equation*}
By Case \ref{case:a>0b=0c>0d=0e=0}) we see that no solution exists.

\item \label{case:a<0b>0c>0d<0e>0} $a<0,b>0,c>0,d<0,e>0\Rightarrow \ga=a=-h, \gb=0, \gc=0, \gd=d=-m,\gg\ge0$. Then
\begin{equation*}
 x^h (x-i)^m -C=D (x-1)^b (x+1)^c(x+i)^e.
\end{equation*}
By Case \ref{case:a>0b=0c=0d>0e=0}) we see that no solution exists.

\item \label{case:a<0b=0c<0d<0e>0} $a<0,b=0,c<0,d<0,e>0\Rightarrow \ga=a=-h, \gb\ge0, \gc=c=-l, \gd=d=-m,\gg\ge0$. Then
\begin{equation*}
 x^h (x+1)^l (x-i)^m-C (x-1)^\gb=D(x+i)^e.
\end{equation*}
Applying $x\to -ix$ in Case \ref{case:a<0b<0c>0d<0e=0}) we see that no solution exists.

\item \label{case:a<0b=0c<0d=0e>0} $a<0,b=0,c<0,d=0,e>0\Rightarrow \ga=a=-h, \gb\ge0, \gc=c=-l, \gd\ge 0,\gg\ge0$. Then
\begin{equation*}
 x^h (x+1)^l-C (x-1)^\gb (x-i)^\gd=D(x+i)^e.
\end{equation*}
Taking complex conjugation of all the coefficients in Case \ref{case:a<0b=0c<0d>0e=0}) we see that no solution exists.

\item \label{case:a<0b=0c=0d<0e>0} $a<0,b=0,c=0,d<0,e>0\Rightarrow \ga=a=-h, \gb\ge0, \gc\ge0, \gd=d=-m,\gg\ge0$. Then
\begin{equation*}
 x^h (x-i)^m-C (x-1)^\gb (x+1)^\gc=D(x+i)^e.
\end{equation*}
Applying $x\to ix$ in Case \ref{case:a<0b>0c<0d=0e=0}) we get $\boxed{f=\frac{i(x+i)}{x(i-x)}}.$

\item \label{case:a<0b=0c<0d>0e>0} $a<0,b=0,c<0,d>0,e>0\Rightarrow \ga=a=-h, \gb\ge0, \gc=c=-l, \gd=0,\gg\ge0$. Then
\begin{equation*}
 x^h (x+1)^l -C (x-1)^\gb=D (x-i)^d(x+i)^e.
\end{equation*}
Applying $x\to -ix$ in Case \ref{case:a<0b>0c>0d<0e=0}) we get $\boxed{f=\frac{x^2+1}{x(x+1)}}.$

\item  \label{case:a<0b=0c>0d<0e>0} $a<0,b=0,c>0,d<0,e>0\Rightarrow \ga=a=-h, \gb\ge0, \gc=0, \gd=d=-m,\gg\ge0$. Then
\begin{equation*}
 x^h (x-i)^m-C (x-1)^\gb=D(x+1)^c(x+i)^e.
\end{equation*}
Applying $x\to -ix$ in Case \ref{case:a<0b<0c>0d>0e=0}) we see that no solution exists.

\end{enumerate}

\bigskip
So far, we have produced 186 4-unital rational functions.

\bigskip
\noindent\  \textbf{XIX}. $e<0$. For every case with this condition we may find the solution
by considering the case with the exponent conditions  $<0$ and $>0$ switched while keeping  "=" unchanged. In this way,
we can produce the other 58 4-unital rational functions with $x=-i$ as a pole.

This concludes the proof of Theorem \ref{thm:U4}.

\section{Proof of Theorem \ref{thm:U3}: complete set of 3-unital functions}
Suppose $f= D x^a (x-1)^b(x-u)^c (x-\baru)^d, g=C x^\ga (x-1)^\gb (x-u)^\gc(x-\baru)^\gd\in\vU_3$ such that $1-g=f$.
We now break into different cases according the signs of the exponents using the more transparent form in \eqref{equ:GenForm}.

\bigskip
\noindent \textbf{I}. $a>0,b,c,d\ge0$.
\begin{enumerate}[{1)}]
\item\label{case:a>0b=0c=0d=0} $a>0,b=c=d=0\Rightarrow \ga=0, \gb,\gc\ge0$. Then
\begin{equation*}
  1-C (x-1)^\gb (x-u)^\gc (x-\baru)^\gd=D x^a.
\end{equation*}

i) If $a=1$, then $\gb+\gc+\gd=10$.

\begin{itemize}
 \item $\gb=1 \Rightarrow \gc=0, \gd=0$. Then $x=1\Rightarrow D=1 \Rightarrow \boxed{f=x}.$

 \item $\gc=1 \Rightarrow \gb=0, \gd=0$. Then $x=u\Rightarrow D=\baru  \Rightarrow\boxed{f=\baru x}.$

 \item $\gd=1 \Rightarrow \gb=0, \gc=0$. Then $x=\baru \Rightarrow D=\baru  \Rightarrow\boxed{f=ux}.$
\end{itemize}

ii) If $a=\gb+\gc+\gd\ge 2\Rightarrow \gb+u\gc+\baru \gd=\gb-\gd+u(\gc-\gd)=0$ by considering 2HD. So $\gb=\gc=\gd\Rightarrow 1-C(x^3-1)^\gb=Dx^a \Rightarrow \gb=1$ since the RHS has only one term. So $C=-1$, $D=1, a=3\Rightarrow \boxed{f=x^3}.$

\item \label{case:a>0b>0c=0d=0} $a>0,b>0,c=d=0\Rightarrow \ga=\gb=0, \gc\ge 0, \gd\ge 0$. Then
\begin{equation*}
  1-C (x-u)^\gc(x-\baru)^\gd=D x^a (x-1)^b.
\end{equation*}
This is impossible by the 3-adic argument.

\item \label{case:a>0b=0c>0d=0} $a>0,b=0,c>0,d=0\Rightarrow \ga=\gc=0, \gb\ge0, \gd\ge0$. Then
\begin{equation*}
  1-C (x-1)^\gb(x-\baru)^\gd=D x^a (x-u)^c.
\end{equation*}
This is impossible by the 3-adic argument.

\item \label{case:a>0b=0c=0d>0} $a>0,b=c=0,d>0$. From the proceeding case this is impossible by symmetry $u\leftrightarrow \baru $.

\item\label{case:a>0b>0c>0d=0} $a>0,b>0,c>0,d=0 \Rightarrow \ga=\gb=\gc=0, \gd\ge 0$.
\begin{equation*}
  1-C (x-\baru)^\gd=D x^a (x-1)^b (x-u)^c.
\end{equation*}
This is impossible by the 3-adic argument.

\item \label{case:a>0b>0c=0d>0} $a>0,b>0,c=0,d>0$. This is impossible by symmetry $u\leftrightarrow \baru $ using the proceeding case.

\item \label{case:a>0b=0c>0d>0} $a>0,b=0,c>0,d>0 \Rightarrow \ga=\gc=\gd=0, \gb\ge 0$.
\begin{equation*}
  1-C (x-1)^\gb =D x^a (x-u)^c(x-\baru)^d.
\end{equation*}
This is impossible by the 3-adic argument.

\item \label{case:a>0b>0c>0d>0} $a>0,b>0,c>0,d>0 \Rightarrow \ga=\gc=\gb=0$. This case is impossible.

\bigskip\noindent\ \hskip-1.1cm \textbf{II}. $a>0,b<0$.

\item \label{case:a>0b<0c=0d=0} $a>0,b<0,c=0,d=0\Rightarrow \ga=0, \gb=b=-k, \gc\ge0, \gd\ge0$. Then
\begin{equation*}
 (x-1)^k-C (x-u)^\gc(x-\baru)^\gd =D x^a .
\end{equation*}

i) If $a=1, k=1$ then $\gc+\gd\le 1$.
\begin{itemize}
 \item $\gc=\gd=0\Rightarrow D=1 \Rightarrow \boxed{f=\frac{x}{x-1}}.$

 \item $\gc=1, \gd=0\Rightarrow (1-C)x-1+Cu=Dx \Rightarrow D=1-\baru  \Rightarrow \boxed{f=\frac{(1-\baru)x}{x-1}}.$

 \item $\gd=1, \gc=0\Rightarrow D=1-u$ by symmetry, $\Rightarrow \boxed{f=\frac{(1-u)x}{x-1}}.$
\end{itemize}

ii) If $a=1, k\ge 2$ then $\gc+\gd=k$ and $C=1$. By considering 2HD we have $-k+u\gc+\baru \gd=-k-\gd+u(\gc-\gd)=D$ (and therefore $k=2$) since $-k-\gd=0$ is impossible. So $\gc=\gd=1$ and
$D=-3 \Rightarrow \boxed{f=\frac{-3x}{(x-1)^2}}.$

iii) If $a\ge 2, k>a$ then $\gc+\gd=k$ and $C=1$. By considering 2HD we have $-k+u\gc+\baru \gd=-k-\gd+u(\gc-\gd)=D$ (and therefore $k=a+1$) since $-k-\gd=0$ is impossible. So $\gc=\gd, k=2\gc\ge 3$ and $D=-3\gc \Rightarrow (x-1)^{2\gc}-(x^2+x+1)^\gc =-3\gc x^{2\gc-1}$. This is impossible by the sign pattern.

iv) If $a\ge 2, k=a$ then $\gc+\gd\le k$

\begin{itemize}
 \item $\gc+\gd= k$.
By considering 2HD we have $-k+C(u\gc+\baru \gd)=0$. On the other hand, $x=0\Rightarrow C=\pm 1, \pm u,\pm \baru $. Thus we must have either $\gc=\gd$, or $\gc=0$, or $\gd=0$. First we assume $C=\pm 1$ and $\gc=\gd$ then the equation becomes $(x-1)^k-C(x^2+x+1)^\gc =D x^k$. So $k=2$ by the sign pattern. This is impossible for $C=\pm 1$. Next we assume $C=\pm u,\gd=0$. Then
the equation becomes $(x-1)^k-\baru (x-u)^k =D x^k$. This is absurd since the coefficient of $x^{k-2}$ on the LHS
is $\binom{k}{2}(1-u^4)\ne 0$. So this case is impossible. By symmetry neither is $\gc=0$.

 \item $\gc+\gd<k \Rightarrow D=1$.
By considering 2HD we have $\gc+\gd=k-1$ and $-k-C=0\Rightarrow C=-k$. Then
by considering 3HD we have $\gc=\gd$ and the equation becomes
$k(x^2+x+1)^\gc =x^k-(x-1)^k$ and $k=2\gc+1\ge 2$. This is absurd by the sign pattern.
\end{itemize}

v) If $a\ge 2, k<a$ then $\gc+\gd=k \Rightarrow C=-D$ by the leading coefficients.
Then $x=1\Rightarrow D(1-u)^\gc(1-\baru)^\gd =D \Rightarrow \gc=\gd=0$ which contradicts the assumption.

\item \label{case:a>0b<0c=0d>0} $a>0,b<0,c=0,d>0\Rightarrow \ga=0, \gb=b=-k, \gc\ge0,\gd=0$. Then
\begin{equation*}
 (x-1)^k-C (x-u)^\gc =D x^a (x-\baru)^d.
\end{equation*}
First $k=\gc$ by the 3-adic argument. Next, $x=0\Rightarrow C=\pm 1, \pm u,\pm \baru $.

i) $C=1$. Then $a+d=k-1\ge2, -k+uk=D$ by 2HD.
By 3HD we get $-\binom{k}{2} \baru =-k(1-u)(-d\baru)\Rightarrow -(k-1)=2d$ which is absurd.

ii) $C\ne 1$. Then $a+d=k, 1-C=D, -k+uCk=-d(1-C)\baru =d(1-C)(u+1)$ by considering 2HD. Then only $C=-u$ is possible $\Rightarrow k=d$. But this implies $a=0$ which is a contradiction.

\item\label{case:a>0b<0c>0d=0} $a>0,b<0,c>0,d=0$. By symmetry $u\leftrightarrow \baru $ we get no solution either.

\item\label{case:a>0b<0c>0d>0} $a>0,b<0,c>0,d>0\Rightarrow \ga=0, \gb=b=-k, \gc=0, \gd=0$. Then
\begin{equation*}
 (x-1)^k-C=D x^a (x-u)^c (x-\baru)^d.
\end{equation*}
This is impossible by the 3-adic argument.

\item \label{case:a>0b<0c<0d=0} $a>0,b<0,c<0,d=0\Rightarrow \ga=0, \gb=b=-k, \gc=c=-l, \gd\ge0$. Then
\begin{equation*}
 (x-1)^k (x-u)^l-C (x-\baru)^\gd=D x^a.
\end{equation*}
First note that $k+l=\gd\ge 2$ by the single power $3$-adic argument. Also, $x=0 \Rightarrow C=\pm 1, \pm u,\pm \baru $.

i) $a=1\Rightarrow C=1$.
\begin{itemize}
 \item $k+l\ge 3$. Then $k+ul- \gd \baru =0$ by 2HD. So $k+\gd+u(l+\gd)=0$ which is absurd.
 \item $k=l=1$. Then $x^2+\baru x+u-(x-\baru)^2=Dx \Rightarrow D=3\baru  \Rightarrow \boxed{f=\frac{3\baru  x}{(x-1)(x-u)}}.$ \end{itemize}

ii) $a\ge 2, k+l=\gd=a \Rightarrow 1-C=D$. Then $k+ul-C\gd \baru =0$ by 2HD.
Clearly $C=\pm u\Rightarrow l=0$ which is impossible. Also $C=\pm \baru \Rightarrow k=0$ which impossible, too.
If $C=-1$, then $k+ul-\gd(u+1)=0 \gd=k=l\Rightarrow 2a=a$ which is absurd.

iii) $a\ge 2, \gd=k+l>a \Rightarrow C=1$. The equation becomes
$$(x-1)^k (x-u)^l-(x-\baru)^{k+l}=D x^a.$$

\begin{itemize}
 \item $a=k+l-1$. Then by 3HD on the LHS
 $\binom{k}{2}+klu+\binom{l}{2}\baru =\binom{a+1}{2}u\Rightarrow \binom{k}{2}-\binom{l}{2}=0, kl-\binom{l}{2}=\binom{a+1}{2}\Rightarrow k=l, a=k$. But $a=k+l-1$ implies $a=1$ which contradicts to $a\ge 2$.

 \item $a<k+l-1$. Then $0=k+ul-(k+l) \baru =2k+l+(2l+k)u$ by 2HD. This is impossible.
\end{itemize}

\item\label{case:a>0b<0c=0d<0} $a>0,b<0,c=0,d<0$. By symmetry $u\leftrightarrow \baru $ we get $\boxed{f=\frac{3ux}{(x-1)(x-\baru  )}}.$

\item\label{case:a>0b<0c<0d>0} $a>0,b<0,c<0,d>0\Rightarrow \ga=0, \gb=b=-k, \gc=c=-l, \gd=0$. Then
\begin{equation*}
 (x-1)^k (x-u)^l-C =D x^a(x-\baru)^d.
\end{equation*}
This is impossible by the 3-adic argument.

\item \label{case:a>0b<0c>0d<0}$a>0,b<0,c>0,d<0$. This is impossible by symmetry $u\leftrightarrow \baru $ using the proceeding case.

\item\label{case:a>0b<0c<0d<0} $a>0,b<0,c<0,d<0\Rightarrow \ga=0, \gb=b=-k, \gc=c=-l, \gd=d=-m$. Then
\begin{equation*}
 (x-1)^k (x-u)^l(x-\baru)^m-C =D x^a.
\end{equation*}
Clearly $k+l+m=a$. By 2HD we get $l=m=k$. Thus the equation becomes
\begin{equation*}
 (x^3-1)^k -C =D x^a.
\end{equation*}
We must have $k=1$. Then $x=0\Rightarrow C=-1; x=1 \Rightarrow D=1$. Thus $\boxed{f=\frac{x^3}{x^3-1}}.$.

\bigskip\noindent\ \hskip-1.1cm \textbf{III}. $a>0,b\ge0,c<0$ or $d<0$.

\item \label{case:a>0b>0c<0d<0} $a>0,b>0,c<0,d<0\Rightarrow \ga=0, \gb=0, \gc=c=-l, \gd=d=-m$. Then
\begin{equation*}
 (x-u)^l(x-\baru)^m-C =D x^a(x-1)^b.
\end{equation*}
This is impossible by the 3-adic argument.

\item \label{case:a>0b>0c<0d=0} $a>0,b>0,c<0,d=0\Rightarrow \ga=0, \gb=0, \gc=c=-l, \gd\ge 0$. Then
\begin{equation*}
 (x-u)^l-C(x-\baru)^\gd =D x^a(x-1)^b.
\end{equation*}
Note that $l=\gd$ by the 3-adic argument. Also, $x=0 \Rightarrow C=\pm 1, \pm u,\pm \baru $.

i) $C=1\Rightarrow lu-l\baru =D$ and $a+b=l-1\ge 2$ by 2HD. By 3HD we have $\binom{l}{2}(\baru -u)=\binom{b}{2}D=\binom{b}{2} l(u-\baru) \Rightarrow \binom{l}{2} =-l\binom{b}{2}$ which is absurd.

ii) $C\ne 1\Rightarrow 1-C=D$ and $l=a+b$. By 2HD we get $lu-l\baru =bD$.
Note $D$ cannot be a unit in $\Z[u]$ by the 3-adic evaluation of $lu-l\baru =bD$,
we see that $C=u$ or $C=\baru $. Thus $lu=b$ or
$lu=b(1+u)=-b\baru $. Neither can hold.

\item \label{case:a>0b>0c=0d<0} $a>0,b>0,c=0,d<0$. No solution by symmetry $u\leftrightarrow \baru $ from the proceeding case.

\item \label{case:a>0b>0c<0d>0} $a>0,b>0,c<0,d>0\Rightarrow \ga=0, \gb=0, \gc=c=-l, \gd=0$. Then
\begin{equation*}
 (x-u)^l-C =D x^a(x-1)^b (x-\baru)^d.
\end{equation*}
This is impossible by the 3-adic argument.

\item \label{case:a>0b>0c>0d<0} $a>0,b>0,c>0,d<0$. No solution by symmetry $u\leftrightarrow \baru $ from the proceeding case.

\item \label{case:a>0b=0c<0d<0} $a>0,b=0,c<0,d<0\Rightarrow \ga=0, \gb\ge0, \gc=c=-l, \gd=d=-m$. Then
\begin{equation*}
 (x-u)^l (x-\baru)^m-C (x-1)^\gb=D x^a.
\end{equation*}
Note that $l+m=\gb$ by the single power 3-adic argument. Also, $x=0 \Rightarrow C=\pm 1, \pm u,\pm \baru $.

i) $l+m=\gb>a\Rightarrow C=1$.
\begin{itemize}
 \item $a= l+m-1$. By 3HD,
 $\binom{\gb}{2} =\binom{l}{2}u+lm+\binom{m}{2}\baru =(\binom{l}{2}-\binom{m}{2})u+lm-\binom{m}{2}\Rightarrow l=m, \gb=m+1$. Thus $\gb=2, l=m=1, a=1.$ By 2HD, $D=\gb-lu-m\baru =3 \Rightarrow \boxed{f=\frac{3x}{x^2+x+1}}.$

 \item $a<l+m-1$. By 2HD, $0=lu+m\baru -\gb=(l-m)u-m-\gb\Rightarrow m+\gb=0$ which is absurd.
\end{itemize}

ii) $l+m=\gb=a$. Then $1-C=D$ so $x=1\Rightarrow (1-u)^l (1-\baru)^m=D$. Then $C=u,\baru $ are the only possible choices. Then the 3-adic value implies that $l+m=1$ which is absurd.

\item \label{case:a>0b=0c<0d=0} $a>0,b=0,c<0,d=0\Rightarrow \ga=0, \gb\ge0, \gc=c=-l, \gd\ge 0$. Then
\begin{equation*}
 (x-u)^l-C (x-1)^\gb (x-\baru)^\gd=D x^a.
\end{equation*}

i) $\gb=\gd=0$. Then we must have $l=a=D=1\Rightarrow \boxed{f=\frac{\baru x}{\baru x-1}}.$

If $\gb+\gd\ne 0$ then $l=\gb+\gd$ by the single power 3-adic argument. Also, $x=0 \Rightarrow C=\pm 1, \pm u,\pm \baru $.

ii) $l=\gb+\gd>a\Rightarrow C=1$.
\begin{itemize}
 \item $a=l-1$. By 3HD,
 $\binom{l}{2}\baru =\binom{\gb}{2}+\gb\gd \baru +\binom{\gd}{2}u \Rightarrow (\binom{l}{2}-\gb\gd+\binom{\gd}{2})\baru =\binom{\gb}{2}-\binom{\gd}{2}\Rightarrow \gb=\gd, l=\gb+1$. Thus $\gb=\gd=1, l=2, a=1$. By 2HD, $D=\gb+\gd \baru -lu=-3u \Rightarrow \boxed{f=\frac{-3ux}{(x-u)^2}}.$

 \item $a<l-1$. By 2HD, $0=\gb+\gd \baru -lu=\gb-\gd-(\gd+l)u\Rightarrow \gd+l=0$ which is absurd.
\end{itemize}

iii) $l=\gb+\gd=a$. Then $1-C=D$ so $x=1\Rightarrow (1-u)^l=D$. Then $C=u,\baru $ are the only possible choices. Then the 3-adic value implies that $l=1\Rightarrow D=1-u$ or $D=1-\baru \Rightarrow \boxed{f=\frac{(1-u)x}{x-u}, \frac{(1-\baru)x}{x-u}}.$

\item \label{case:a>0b=0c=0d<0} $a>0,b=0,c=0,d<0$. By symmetry $u\leftrightarrow \baru $ we get
\begin{equation*}
\boxed{f=\frac{ux}{ux-1},\quad \frac{-3\baru  x}{(x-\baru)^2},\quad \frac{(1-\baru)x}{x-\baru},\quad  \frac{(1-u)x}{x-\baru}}.
\end{equation*}

\item \label{case:a>0b=0c<0d>0} $a>0,b=0,c<0,d>0\Rightarrow \ga=0, \gb\ge0, \gc=c=-l, \gd=0$. Then
\begin{equation*}
 (x-u)^l -C (x-1)^\gb =D x^a (x-\baru)^d.
\end{equation*}

First, $l=\gb$ by the 3-adic argument. Thus $l=\gb\ge a+d\ge 2$.

i) $l=\gb>a+d\Rightarrow C=1$.
\begin{itemize}
 \item $a+d=l-1$. By 2HD, $l-lu=D$.
 By 3HD,
 $\binom{l}{2}\baru -\binom{l}{2}=-Dd\baru =(lu-l)d\baru =ld(1-\baru)\Rightarrow (\binom{l}{2}+ld)(\baru -1)=0$, which is absurd.

 \item $a+d<l-1$. By 2HD, $l-lu=0$ which is absurd.
\end{itemize}

ii) $l=\gb=a+d\ge 2$. Then $1-C=D$ so $x=1\Rightarrow (1-u)^l =D(1-\baru)^d$. Since $d<l$ 3-adic evaluation of $D$ is $<1$ so that $C=u,\baru $ are the only possible choices, which implies that $a=1$. The equation becomes
$(x-u)^l-C (x-1)^l =(1-C)x(x-\baru)^{l-1}$. Then second highest degree term yields $l(C-u)=(C-1)(l-1) \baru $, which
is satisfied by neither $C=u$ $(\Rightarrow l=1)$ nor $C=\baru $ ($\Rightarrow l=1-l$).

\item \label{case:a>0b=0c>0d<0} $a>0,b=0,c>0,d<0$. No solution by symmetry $u\leftrightarrow \baru $ from the proceeding case.

\bigskip\noindent\ \hskip-1.1cm \textbf{IV}. $a=0,b,c,d\ge0$. We set $1/D=-B$. Then we must have $C=-D$.

\item \label{case:a=0b=0c=0d=0} $a=0,b=c=d=0\Rightarrow \ga\ge 0, \gb,\gc\ge0$. This is clearly impossible.

\item \label{case:a=0b>0c=0d=0} $a=0,b>0, c=d=0\Rightarrow \ga=\gb=0, \gc\ge 0, \gd\ge 0$. Then we must have
\begin{equation*}
  B-x^\ga (x-u)^\gc(x-\baru)^\gd=-(x-1)^b.
\end{equation*}

i) $\ga\ge 1$. This is impossible by the 3-adic argument except for $\gc=\gd=0$ in which case we must have
$b=1$ by the sign pattern. Thus $B=\ga=1 \Rightarrow \boxed{f=1-x}.$
ii) $\ga=0$.
\begin{itemize}
 \item $\gc=1,\gd=0$. Then $b=1$ and $B=1-u \Rightarrow \boxed{f=\frac{x-1}{u-1}}.$ \item $\gc=0,\gd=1$. Then $b=1$ and $B=1-\baru  \Rightarrow \boxed{f=\frac{x-1}{\baru -1}}.$

 \item $\gc+\gd\ge 2$. Then by 2HD $\gc u+\gd \baru =b\Rightarrow (\gc-\gd)u=b+\gd=0$ which is absurd.
\end{itemize}

\item \label{case:a=0b=0c>0d=0} $a=0,b=0,c>0, d=0\Rightarrow \ga=\gc=0, \gb\ge0, \gd\ge0$. Then we must have
\begin{equation*}
  B-x^\ga (x-1)^\gb(x-\baru)^\gd=-(x-u)^c.
\end{equation*}

i) $\ga\ge 1$. This is impossible by the 3-adic argument except for $\gb=\gd=0$ in which case we must have
$c=1$ the LHS has only two terms. Thus $B=u,\ga=1 \Rightarrow \boxed{f=1-\baru x}.$

ii) $\ga=0$.
\begin{itemize}
 \item $\gb=0,\gd=1$. Then $c=1$ and $B=u-\baru  \Rightarrow \boxed{f=\frac{x-u}{\baru -u}}.$

 \item $\gb=1,\gd=0$. Then $c=1$ and $B=u-1 \Rightarrow \boxed{f=\frac{x-u}{1-u}}.$

 \item $\gc+\gd\ge 2$. Then by 2HD $\gc+\gd \baru =c u\Rightarrow \gc-\gd=(c+\gd)u=0$ which is absurd.
\end{itemize}

\item \label{case:a=0b=0c=0d>0} $a=0,b=c=0,d>0$. By symmetry $u\leftrightarrow \baru $ from above we get
\begin{equation*}
 \boxed{f=1-ux, \qquad \frac{x-\baru }{u-\baru }, \qquad \frac{x-\baru }{1-\baru}}.
\end{equation*}

\item \label{case:a=0b>0c>0d=0} $a=0,b>0, c>0, d=0 \Rightarrow \ga=\gb=\gc=0, \gd\ge 0$. Then we must have
\begin{equation*}
  B-x^\ga (x-\baru)^\gd=-(x-1)^b (x-u)^c.
\end{equation*}
i) $\ga\ge 1$. This is impossible by the 3-adic argument except for $\gd=0$. But the RHS has at least three terms,
i.e., $x^{b+c}, -(b+cu)x^{b+c-1}$ and the nonzero constant term,
while the LHS has only two terms, which is a contradiction.

ii) $\ga=0$. Then $\gd=b+c\ge 2$ and by 2HD $\gd \baru =b+cu\Rightarrow b+\gd=-(c+\gd)u=0$,
which is impossible.

\item \label{case:a=0b>0c=0d>0} $a=0,b>0,c=0, d>0$. No solution by symmetry $u\leftrightarrow \baru $ using the proceeding case.

\item \label{case:a=0b=0c>0d>0} $a=0,b=0,c>0, d>0 \Rightarrow \ga=\gc=\gd=0, \gb\ge 0$. Then we must have
\begin{equation*}
  B-x^\ga (x-1)^\gb =- (x-u)^c(x-\baru)^d.
\end{equation*}
i) $\ga\ge 1$. This is impossible by the 3-adic argument except for $\gb=0$. But the RHS has at least three terms,
i.e., $x^{c+d}, -(d\baru +cu)x^{c+d-1}$ and the nonzero constant term,
while the LHS has only two terms, which is a contradiction.

ii) $\ga=0$. $\gb=c+d\ge 2$ and by 2HD $\gb=cu+d\baru \Rightarrow \gb+d=(c-d)u=0$,
which is impossible.

\item \label{case:a=0b>0c>0d>0} $a=0,b>0,c>0,d>0 \Rightarrow \ga=\gc=\gb=0$. Then we must have
\begin{equation*}
  B- x^\ga =-(x-1)^b (x-u)^c (x-\baru)^d.
\end{equation*}
$x=1\Rightarrow B=1$. Then $x=u\Rightarrow 3|\ga=3n$. So the LHS is invariant under $x\to xu$ and $x\to x\baru $.
Thus we have $b=c=d$ so that $1-x^{3n}=-(x^3-1)^b$. Thus $b=n=1$ and we get $\boxed{f=1-x^3}.$

\newpage
\noindent\ \hskip-1.1cm \textbf{V}. $a=0,b<0$.

\item \label{case:a=0b<0c=0d=0} $a=0,b<0,c=0,d=0\Rightarrow \ga\ge 0, \gb=b=-k, \gc\ge0, \gd\ge0$. Then
\begin{equation*}
 (x-1)^k-C x^\ga (x-u)^\gc(x-\baru)^\gd =D .
\end{equation*}
By Case \ref{case:a=0b>0c=0d=0}) we get
\begin{equation*}
 \boxed{f=\frac{1}{1-x}, \qquad \frac{u-1}{x-1},\qquad \frac{\baru -1}{x-1}}.
\end{equation*}

\item \label{case:a=0b<0c=0d>0} $a=0,b<0,c=0,d>0\Rightarrow \ga\ge 0, \gb=b=-k, \gc\ge0,\gd=0$. Then
\begin{equation*}
 (x-1)^k-C x^\ga (x-u)^\gc =D (x-\baru)^d.
\end{equation*}

i) $\ga=\gc=0$. Then $D=1$ and $k=d$. If $k=d\ge 2$ then by 2HD $k=d\baru $ which is absurd.
So $k=d=1$ and $\boxed{f=\frac{x-\baru }{x-1}}.$

ii) $\ga\ge 1, \ga+\gc=k$. Then $1-C=D$ and $k=d$ by the single power 3-adic argument. Also
$x=0\Rightarrow D=\pm 1, \pm u,\pm \baru $. By 2HD
$-k+(1-D)\gc u=-Dd \baru  \Rightarrow (\gc-\gd)u=k+\gd=0$.

\begin{itemize}
 \item $D=-1$. Then $-k+2\gc u=d \baru \Rightarrow (2\gc+d)u=k-d=0$ which is absurd.

 \item $D=\pm u$. Then $\gc u\pm \gc(u+1)\pm k=k \Rightarrow -\gc-k=k$ or $\gc(2u+1)=0$. So $D=-u$ is the only possible choice and $\gc=0$. If $k\ge 2$ then by 3HD $\binom{k}{2}=-u \binom{k}{2} u^4$
 which is absurd. Thus $k=1 \Rightarrow \boxed{f=\frac{1-ux}{1-x}}.$.

 \item $D=\pm \baru $. Then $\gc u \pm k u=k\pm \gc$. So $D=-\baru $ is the only possible choice and $\gc=k$.
 If $k\ge 2$ then by 3HD $\binom{k}{2}=-\baru  \binom{k}{2} u^4$ which is absurd. Thus
 $k=1 \Rightarrow \boxed{f=\frac{u-\baru x}{x-1}}.$.
\end{itemize}

iii) $\ga\ge 1, \ga+\gc<k$. Then $d=k$ and $D=1$.
\begin{itemize}
 \item $\ga+\gc=k-1$. Then by 2HD $k+C=k \baru  \Rightarrow C=k(\baru -1)$. If $k\ge 2$ then by 3HD $\binom{k}{2} u^4=\binom{k}{2}+C\gc u=\binom{k}{2}+k\gc(1-u)\Rightarrow \binom{k}{2}+k\gc=0$
 which is absurd. Thus $k=d=1 \Rightarrow \ga=0$ which contradicts the assumption.

 \item $\ga+\gc<k-1$. Then by 2HD $k=k \baru $ which is impossible.
\end{itemize}

iv) $\ga=0, \gc\ge 1$.
\begin{itemize}
 \item $k=\gc=d$. Then $1-C=D$ and by 2HD $-k+(1-D)k u=-D k \baru =Dk(u+1) \Rightarrow (1-2D)u=D+1=0$, which has no solution.

 \item $k=\gc>d$. Then $C=1$. If $d<k-1$ then by 2HD $-k+k u=0$ which is absurd.
 So $d=k-1$ and $-k+k u=D.$ But if further $k\ge 2$ then by 3HD
 $$\binom{k}{2}(1-u)=-D(k-1)\baru =(ku-k)(k-1)(u+1) =k(k-1)(\baru -1)=k(k-1)(-u-2)$$
 which is impossible.

 \item $k=d>\gc$. Then $D=1$. No solution from the subcase immediately above by the symmetry $u\leftrightarrow \baru $.

 \item $k<d=\gc$. Then $d=\gc$ and $C=-D$. Let $B=1/D$. The equation becomes $B(x-1)^k=(x-\baru)^d-(x-u)^d$.
 If $k<d-1$ then by 2HD $d(u-\baru)\ne 0$ on the RHS which is impossible.
 If $k=d-1\ge 1$ then $B=d(u-\baru)$ and $-Bk=\binom{d}{2}(u^4-\baru)$ by the second and third highest degree term,
 which implies that $-dk=\binom{d}{2}$, again impossible.

\end{itemize}

\item \label{case:a=0b<0c>0d=0} $a=0,b<0, c>0,d=0$. By symmetry $u\leftrightarrow \baru $ there are three possible solutions:
\begin{equation*}
 \boxed{f=\frac{x-u}{x-1}, \qquad \frac{1-\baru x}{1-x}, \qquad \frac{\baru -ux}{x-1}}.
\end{equation*}

\item \label{case:a=0b<0c>0d>0} $a=0,b<0, c>0,d>0\Rightarrow \ga\ge 0, \gb=b=-k, \gc=0, \gd=0$. Then
\begin{equation*}
 (x-1)^k-Cx^\ga=D (x-u)^c (x-\baru)^d.
\end{equation*}

i) $\ga\ge 1$. By the single power $3$-adic argument, $k=c+d\ge2$ and $D=\pm 1, \pm u,\pm \baru $.
\begin{itemize}
 \item $\ga=k$. Then $1-C=D$ and 2HD and 3HD $\Rightarrow k=D(cu+d\baru)$ which cannot hold for any of the possible values of $D$.
 \item $\ga=k-1$. Then $D=1$ and 3HD $\Rightarrow \binom{k}{2}=\big(\binom{c}{2}\baru +cd+\binom{d}{2}u\big)\Rightarrow \binom{k}{2}=cd-\binom{c}{2}, \binom{c}{2}=\binom{d}{2}\Rightarrow c=d,k=c+1\Rightarrow c=d=1, k=2 \Rightarrow \boxed{f=\frac{x^2+x+1}{(x-1)^2}}.$
 \item $\ga<k-1$. Then $D=1$ and 2HD $\Rightarrow k=(cu+d\baru)=u(c-d)-d\Rightarrow k=-d$ which is impossible.
\end{itemize}

ii) $\ga=0$. It's clear that $D=1$ and $k=c+d\ge 2$. By 2HD $k=D(cu+d\baru)$. By 3HD $\binom{k}{2}=D\big(\binom{c}{2}\baru +\binom{d}{2}u\big)$. Thus $(c+du)\binom{k}{2}=k\big(\binom{c}{2}u+\binom{d}{2}\big)$, or $(c+du)(k-1)=c(c-1)u+d(d-1)\Rightarrow c(k-1)=d(d-1), d(k-1)=c(c-1)\Rightarrow c=d=k$ contradicting to $k=c+d$.

\item \label{case:a=0b<0c<0d=0} $a=0,b<0, c<0,d=0\Rightarrow \ga\ge 0, \gb=b=-k, \gc=c=-l, \gd\ge0$. Then
\begin{equation*}
 (x-1)^k (x-u)^l-C x^\ga (x-\baru)^\gd=D.
\end{equation*}
Clearly $C=1$. Observe that $\ga\gd=0$ by the 3-adic argument.

i) $\ga=0$. By 2HD $k+lu-\gd \baru =k+lu+\gd(u+1)=0\Rightarrow k+\gd=0$ which is absurd.

ii) $\gd=0$. By 2HD $k+lu=0$ which is impossible either.

\item $a=0,b<0, c=0,d<0$. By symmetry $u\leftrightarrow \baru $ we get no solution.

\item \label{case:a=0b<0c<0d>0} $a=0,b<0, c<0,d>0\Rightarrow \ga\ge 0, \gb=b=-k, \gc=c=-l, \gd=0$. Then
\begin{equation*}
 (x-1)^k (x-u)^l-C x^\ga =D(x-\baru)^d.
\end{equation*}

i) If $\ga=0$ then $D=1$ and by 2HD $k+lu=d\baru =-d(u+1)$ which is impossible.

ii) If $\ga\ge 1$ then $k+l=d$ by the single power 3-adic argument. Also $x=0\Rightarrow D=\pm 1, \pm u,\pm \baru $.
\begin{itemize}
 \item $\ga=d$. Then $1-C=D$. By 2HD $k+lu=Dd\baru =D(l+d)\baru $ which cannot hold for any of the possible values of D.

 \item $\ga<d$. Then $D=1$ and by
 3HD $\binom{k}2+klu+\binom{l}2 \baru =\binom{d}2 u\Rightarrow \binom{k}2=\binom{l}2, kl-\binom{l}2=\binom{d}2\Rightarrow k=l, d=l+1$. Thus $k=l=1, d=2 \Rightarrow \boxed{f=\frac{(x-\baru)^2}{(x-1)(x-u)}}.$
\end{itemize}

\item \label{case:a=0b<0c>0d<0} $a=0,b<0, c>0,d<0$. By symmetry $u\leftrightarrow \baru $ we get $\boxed{f=\frac{(x-u)^2}{(x-1)(x-\baru)}}.$

\item \label{case:a=0b<0c<0d<0} $a=0,b<0,c<0,d<0\Rightarrow \ga\ge 0, \gb=b=-k, \gc=c=-l, \gd=d=-m$. Then
\begin{equation*}
 (x-1)^k (x-u)^l(x-\baru)^m-C x^\ga =D.
\end{equation*}
By Case \ref{case:a=0b>0c>0d>0}) we immediately get $k=l=m=C=1, D=-1$ and $\boxed{f=\frac{1}{1-x^3}}.$

\bigskip\noindent\ \hskip-1.1cm \textbf{VI}. $a=0,b\ge 0,c<0$ or $d<0$.

\item \label{case:a=0b>0c<0d<0} $a=0,b>0, c<0,d<0\Rightarrow \ga\ge 0, \gb=0, \gc=c=-l, \gd=d=-m$. Then
\begin{equation*}
 (x-u)^l(x-\baru)^m-C x^\ga =D(x-1)^b.
\end{equation*}
By Case \ref{case:a=0b<0c>0d>0}) we have $\boxed{f=\frac{(x-1)^2}{x^2+x+1}}.$

\item $a=0,b>0,c<0,d=0\Rightarrow \ga\ge 0, \gb=0, \gc=c=-l, \gd\ge 0$. Then
\begin{equation*}
 (x-u)^l-C(x-\baru)^\gd =D(x-1)^b.
\end{equation*}
By Case \ref{case:a=0b<0c>0d=0}) we have
\begin{equation*}
 \boxed{f=\frac{x-1}{x-u}, \qquad \frac{1-x}{1-\baru x}, \qquad \frac{x-1}{\baru -ux}}.
\end{equation*}

\item $a=0,b>0,c=0,d<0$. By symmetry $u\leftrightarrow \baru $ we get
\begin{equation*}
 \boxed{f=\frac{x-1}{x-\baru}, \qquad \frac{1-x}{1-ux}, \qquad \frac{x-1}{u-\baru x}}.
\end{equation*}

\item $a=0,b>0,c<0,d>0\Rightarrow \ga\ge 0, \gb=0, \gc=c=-l, \gd=0$. Then
\begin{equation*}
 (x-u)^l-C x^\ga =D(x-1)^b (x-\baru)^d.
\end{equation*}
By Case \ref{case:a=0b<0c>0d<0}) we have $\boxed{f=\frac{(x-1)(x-\baru)}{(x-u)^2}}.$

\item $a=0,b>0,c>0,d<0$. By symmetry $u\leftrightarrow \baru $ we get $\boxed{f=\frac{(x-1)(x-u)}{(x-\baru)^2}}.$

\item $a=0,b=0,c<0,d<0\Rightarrow \ga\ge 0, \gb\ge0, \gc=c=-l, \gd=d=-m$. Then
\begin{equation*}
 (x-u)^l (x-\baru)^m-C x^\ga (x-1)^\gb=D.
\end{equation*}
By Case \ref{case:a=0b=0c>0d>0}) there is no solution.

\item $a=0,b=0,c<0,d=0\Rightarrow \ga\ge 0, \gb\ge0, \gc=c=-l, \gd\ge 0$. Then
\begin{equation*}
 (x-u)^l-C x^\ga (x-1)^\gb (x-\baru)^\gd=D.
\end{equation*}
By Case \ref{case:a=0b=0c>0d=0}) we have
\begin{equation*}
\boxed{f=\frac1{1-\baru x},\qquad \frac{\baru -u}{x-u},\qquad \frac{1-u}{x-u}}.
\end{equation*}

\item $a=0,b=0,c=0,d<0$. By symmetry $u\leftrightarrow \baru $ we get
\begin{equation*}
\boxed{f=\frac{1}{1-ux},\qquad \frac{u-\baru}{x-\baru},\qquad \frac{1-\baru }{x-\baru }}.
\end{equation*}

\item $a=0,b=0,c<0,d>0\Rightarrow \ga\ge 0, \gb\ge0, \gc=c=-l, \gd=0$. Then
\begin{equation*}
 (x-u)^l -C x^\ga (x-1)^\gb =D (x-\baru)^d.
\end{equation*}

i) $\ga=\gb=0$. The $D=1, l=d$ and by 2HD we get $l=d=1$. Thus $\boxed{f=\frac{x-\baru }{x-u}}.$

ii) $\ga\ge1$. Note that $l=d$ by single power 3-adic argument.
\begin{itemize}
 \item $\ga+\gb=l=1$. Then $1-C=D$. By 2HD $D\baru =u-C\gb=u+(1-D)\gb$. If $\gb=0$ then $D=\baru \Rightarrow\boxed{f=\frac{x-\baru }{u(x-u)}}$. If $\gb=1$ then $D=-u\Rightarrow\boxed{f=\frac{1-ux}{x-u}}$.

 \item $\ga+\gb=l\ge 2$. Then $1-C=D$. By 2HD and 3HD $du=Dd\baru , \binom{d}{2}=D\binom{d}{2}u$ which is impossible.

 \item $\ga+\gb<l\ge 2$. Then $D=1$ and by 2HD $du=d\baru $, which is absurd.
\end{itemize}

\item $a=0,b=0,c>0,d<0$. By symmetry $u\leftrightarrow \baru $ we get
\begin{equation*}
 \boxed{f=\frac{x-u}{x-\baru },\qquad \frac{u(x-u)}{x-\baru},\qquad \frac{x-u}{1-ux}}.
\end{equation*}

\bigskip\noindent\ \hskip-1.1cm \textbf{VII}. $a<0,b,c,d\ge0$.

\item $a<0,b=c=d=0\Rightarrow \ga=a=-h, \gb,\gc\ge0$. Then
\begin{equation*}
  x^h-C (x-1)^\gb (x-u)^\gc (x-\baru)^\gd=D.
\end{equation*}
By Case \ref{case:a>0b=0c=0d=0}) we have four solutions
\begin{equation*}
\boxed{f=\frac1x,\qquad \frac1{\baru x},\qquad \frac1{ux},\qquad \frac1{x^3}}.
\end{equation*}

\item $a<0,b>0,c=d=0\Rightarrow \ga=\gb=0, \gc\ge 0, \gd\ge 0$. Then
\begin{equation*}
  x^h-C (x-u)^\gc(x-\baru)^\gd=D (x-1)^b.
\end{equation*}
By Case \ref{case:a>0b<0c=0d=0}) we have four solutions
\begin{equation*}
\boxed{f=\frac{x-1}{x},\qquad \frac{x-1}{(1-\baru)x},\qquad \frac{x-1}{(1-u)x},\qquad \frac{(x-1)^2}{-3x}}.
\end{equation*}

\item $a<0,b=0,c>0,d=0\Rightarrow \ga=\gc=0, \gb\ge0, \gd\ge0$. Then
\begin{equation*}
  x^h-C (x-1)^\gb(x-\baru)^\gd=D (x-u)^c.
\end{equation*}
By Case \ref{case:a>0b=0c<0d=0}) we have
\begin{equation*}
\boxed{f=\frac{\baru x-1}{\baru x}, \qquad \frac{(x-u)^2}{-3ux},
\qquad \frac{x-u}{(1-u)x},  \frac{x-u}{(1-\baru)x}}.
\end{equation*}

\item $a=<0,b=c=0, d>0$. By symmetry $u\leftrightarrow \baru $ we have
\begin{equation*}
\boxed{f=\frac{ux-1}{ux}, \qquad \frac{(x-\baru)^2}{-3\baru x}, \qquad \frac{x-\baru }{(1-\baru)x},  \frac{x-\baru}{(1-u)x}}.
\end{equation*}

\item $a<0,b>0,c>0,d=0 \Rightarrow \ga=\gb=\gc=0, \gd\ge 0$.
\begin{equation*}
  x^h-C (x-\baru)^\gd=D (x-1)^b (x-u)^c.
\end{equation*}
By Case \ref{case:a>0b<0c<0d=0}) we have $\boxed{f=\frac{(x-1)(x-u)}{3\baru  x}}.$

\item $a<0,b>0,c=0,d>0$. By symmetry $u\leftrightarrow \baru $ we get $\boxed{f=\frac{(x-1)(x-\baru)}{3u x}}.$

\item $a<0,b=0,c>0,d>0 \Rightarrow \ga=\gc=\gd=0, \gb\ge 0$.
\begin{equation*}
  x^h-C (x-1)^\gb =D (x-u)^c(x-\baru)^d.
\end{equation*}
By Case \ref{case:a>0b=0c<0d<0}) we get $\boxed{f=\frac{x^2+x+1}{3x}}$.

\item $a<0,b>0,c>0,d>0 \Rightarrow \ga=\gc=\gb=0$.
\begin{equation*}
  x^h-C =D (x-1)^b (x-u)^c(x-\baru)^d.
\end{equation*}
By Case \ref{case:a>0b<0c<0d<0}) we get $\boxed{f=\frac{x^3-1}{x^3}}$.

\bigskip\noindent\ \hskip-1.1cm \textbf{VIII}. $a<0,b<0$.

\item $a<0,b<0, c=0,d=0\Rightarrow \ga=a=-h, \gb=b=-k, \gc\ge0, \gd\ge0$. Then
\begin{equation*}
 x^h (x-1)^k-C (x-u)^\gc(x-\baru)^\gd =D .
\end{equation*}
By Case \ref{case:a>0b>0c=0d=0}) there is no solution.

\item $a<0,b<0,c=0,d>0\Rightarrow \ga=a=-h, \gb=b=-k, \gc\ge0,\gd=0$. Then
\begin{equation*}
 x^h (x-1)^k-C (x-u)^\gc =D (x-\baru)^d.
\end{equation*}
By Case \ref{case:a>0b>0c=0d<0}) there is no solution.

\item $a<0,b<0, c>0,d=0$. No solution by symmetry $u\leftrightarrow \baru $ from the proceeding case.

\item $a<0,b<0, c>0,d>0\Rightarrow \ga=a=-h, \gb=b=-k, \gc=0, \gd=0$. Then
\begin{equation*}
 x^h (x-1)^k-C=D (x-u)^c (x-\baru)^d.
\end{equation*}
By Case \ref{case:a>0b>0c<0d<0}) there is no solution.

\item $a<0,b<0, c<0,d=0\Rightarrow \ga=a=-h, \gb=b=-k, \gc=c=-l, \gd\ge0$. Then
\begin{equation*}
 x^h (x-1)^k (x-u)^l-C (x-\baru)^\gd=D.
\end{equation*}
By Case \ref{case:a>0b>0c>0d=0}) there is no solution.

\item $a<0,b<0, c=0,d<0$. No solution by symmetry $u\leftrightarrow \baru $ from the proceeding case.

\item $a<0,b<0, c<0,d>0\Rightarrow \ga=a=-h, \gb=b=-k, \gc=c=-l, \gd=0$. Then
\begin{equation*}
 x^h (x-1)^k (x-u)^l-C =D(x-\baru)^d.
\end{equation*}
By Case \ref{case:a>0b>0c>0d<0}) there is no solution.

\item $a<0,b<0, c>0,d<0$. No solution by symmetry $u\leftrightarrow \baru $ from the proceeding case.

\item $a<0,b<0, c<0,d<0$. This is impossible.

\bigskip\noindent\ \hskip-1.1cm \textbf{IX}. $a<0,b\ge 0,c<0$ or $d<0$.

\item $a<0,b>0, c<0,d<0\Rightarrow \ga=a=-h, \gb=0, \gc=c=-l, \gd=d=-m$. Then
\begin{equation*}
 x^h (x-u)^l(x-\baru)^m-C =D(x-1)^b.
\end{equation*}
By Case \ref{case:a>0b<0c>0d>0}) there is no solution.

\item $a<0,b>0, c<0,d=0\Rightarrow \ga=a=-h, \gb=0, \gc=c=-l, \gd\ge 0$. Then
\begin{equation*}
 x^h (x-u)^l-C(x-\baru)^\gd =D(x-1)^b.
\end{equation*}
By Case \ref{case:a>0b<0c>0d=0}) there is no solution.

\item $a<0,b>0, c=0,d<0$. No solution by symmetry $u\leftrightarrow \baru $ from the proceeding case.

\item $a<0,b>0, c<0,d>0\Rightarrow \ga=a=-h, \gb=0, \gc=c=-l, \gd=0$. Then
\begin{equation*}
 x^h (x-u)^l-C =D(x-1)^b (x-\baru)^d.
\end{equation*}
By Case \ref{case:a>0b<0c>0d<0}) there is no solution.

\item $a<0,b>0, c>0,d<0$. No solution by symmetry $u\leftrightarrow \baru $ from the proceeding case.

\item $a<0,b=0,c<0,d<0\Rightarrow \ga=a=-h, \gb\ge0, \gc=c=-l, \gd=d=-m$. Then
\begin{equation*}
 x^h (x-u)^l (x-\baru)^m-C (x-1)^\gb=D.
\end{equation*}
By Case \ref{case:a>0b=0c>0d>0}) there is no solution.

\item $a<0,b=0,c<0,d=0\Rightarrow \ga=a=-h, \gb\ge0, \gc=c=-l, \gd\ge 0$. Then
\begin{equation*}
 x^h (x-u)^l-C (x-1)^\gb (x-\baru)^\gd=D.
\end{equation*}
By Case \ref{case:a>0b=0c>0d=0}) there is no solution.

\item $a<0,b=0,c=0,d<0$. No solution by symmetry $u\leftrightarrow \baru $ from the proceeding case.

\item $a<0,b=0,c<0,d>0\Rightarrow \ga=a=-h, \gb\ge0, \gc=c=-l, \gd=0$. Then
\begin{equation*}
 x^h (x-u)^l -C (x-1)^\gb =D (x-\baru)^d.
\end{equation*}
By Case \ref{case:a>0b=0c>0d<0}) there is no solution.

\item $a<0,b=0,c>0,d<0$. No solution by symmetry $u\leftrightarrow \baru $ from the proceeding case.

\end{enumerate}

This concludes the proof of Theorem \ref{thm:U3}.

\section*{Appendix A. All $1$-, $2$- and $4$-unital functions}
It is easy to see that the only $1$-unital are given by
$$\vU_1=\langle x\rangle_6=\left\{x, \frac{1}{x}, \frac{x}{x-1}, \frac{x-1}{x}, 1-x, \frac{1}{1-x}\right\}.
$$
In this paper we have proved that there are 36 $2$-unital functions which are given by
\begin{align*}
\vU_2&=\langle x\rangle_6\cup \langle -x\rangle_6\cup \langle x^2\rangle_6\cup \left\langle \frac{1+x}{2} \right\rangle_{\hskip-4pt 6}
\cup \left\langle \frac{1+x}{2x}\right\rangle_{\hskip-4pt 6}  \cup \left\langle \frac{(1+x)^2}{4x}\right\rangle_{\hskip-4pt 6} . \\
&= \left\{
x, \frac{1}{x}, \frac{x}{x-1}, \frac{x-1}{x}, 1-x, \frac{1}{1-x}, \frac{1}{x^2}, x^2, \frac{x^2}{x^2-1}, \frac{x^2-1}{x^2}, \frac{1}{1-x^2},  1-x^2, \right.\\
&\frac{1}{1+x}, \frac{x}{1+x}, \frac{1+x}{x}, -x, \frac{-1}{x}, 1+x, \frac{x+1}2,\frac{1-x}2, \frac{1 + x}{x - 1}, \frac{x - 1}{1 + x}, \frac{2}{1 + x}, \frac{2}{1-x}\\
&\left.\frac{1+x}{1-x}, \frac{1-x}{1+x}, \frac{2x}{1+x}, \frac{2x}{x-1}, \frac{1+ x}{2x}, \frac{x-1}{2x},
\frac{(1+x)^2}{(x-1)^2}, \frac{(x-1)^2}{(1+x)^2}, \frac{(1+x)^2}{4x}, \frac{4x}{(1+x)^2}, \frac{-4x}{(x-1)^2}, \frac{(x-1)^2}{-4x} . \right\}
\end{align*}
This is a corollary of the following complete list of 252 $4$-unital functions
\begin{align*}
\vU_4:= O(x) \cup O(x^2) \cup O(x^4) \cup O\left(\frac{2x}{x^2+1}\right)\cup \left\langle \frac{4x^2}{(x^2+1)^2}\right\rangle_{\hskip-4pt 6} \cup O\left(\frac{2x}{x+1}\right) \cup O\left(\frac{x(x-1)}{x^2+1}\right),
\end{align*}
where by putting $\om=1-i$
\begin{align*}
\ \hskip-2cm
O(x):=& \bigcup_{\eps =\pm 1} \langle \eps x\rangle_6\cup
\langle \eps ix\rangle_6\cup
\left\langle \frac{i(x+1)}{\eps(x-1)}\right\rangle_{\hskip-4pt 6}\cup
\left\langle \frac{i(x-i)}{\eps(x+i)}\right\rangle_{\hskip-4pt 6}\cup
\left\langle \frac{x+\eps}{x-i}\right\rangle_{\hskip-4pt 6}\cup
\left\langle \frac{x+\eps}{x+i}\right\rangle_{\hskip-4pt 6}\\
=
&\left\{ 1-x,x,\frac{1}{1-x},\frac{1}{x},\frac{x-1}{x},\frac{x}{x-1},
\frac{i(x+1)}{x-1},\frac{x-1}{i(x+1)},\frac{\bom(x-i)}{x+1},\frac{\om(x+1)}{2(x-i)},\frac{\om(x-i)}{x-1},\frac{\bom(x-1)}{2(x-i)},\right.\\
&-x,x+1,\frac{-1}{x},\frac{1}{x+1},\frac{x+1}{x},\frac{x}{x+1},
\frac{1-x}{i(x+1)},\frac{i(x+1)}{1-x},\frac{\bom(x+1)}{2(x+i)},\frac{\om(x+i)}{x+1},\frac{\om(x-1)}{2(x+i)},\frac{\bom(x+i)}{x-1},\\
&ix,\frac{-i}{x},\frac{i}{x+i},\frac{x+i}{i},\frac{x+i}{x},\frac{x}{x+i},
\frac{i(x-i)}{x+i},\frac{x+i}{i(x-i)},\frac{\om(x-1)}{x+i},\frac{\bom(x+i)}{2(x-1)},\frac{\bom(x-1)}{x-i},\frac{\om(x-i)}{2(x-1)},\\
&-ix,\frac{-i}{x-i},\frac{i}{x},\frac{x-i}{-i},\frac{x-i}{x},\frac{x}{x-i},
\frac{x-i}{i(x+i)},\frac{i(x+i)}{x-i},\frac{\bom(x-i)}{2(x+1)},\frac{\om(x+1)}{x-i},\frac{\om(x+i)}{2(x+1)},\frac{\bom(x+1)}{x+i},\\
&\frac{x-1}{x-i},\frac{x-i}{x-1},\frac{\om}{1-x},\frac{1-x}{\om},\frac{\om}{x-i},\frac{x-i}{\om},
\frac{x+1}{x-i},\frac{x-i}{x+1},\frac{i-x}{\bom},\frac{\bom}{i-x},\frac{x+1}{\bom},\frac{\bom}{x+1},\\
&\frac{x-1}{x+i},\frac{x+i}{x-1},\frac{\bom}{1-x},\frac{1-x}{\bom},\frac{\bom}{x+i},\frac{x+i}{\bom}
\left.\frac{x+1}{x+i},\frac{x+i}{x+1},\frac{x+1}{\om},\frac{\om}{x+1},\frac{x+i}{-\om},\frac{-\om}{x+i}\right\},\\
\ \hskip-2cm
 O(x^2):=&\bigcup_{\eps =\pm 1}
\langle \eps x^2\rangle_6\cup
\left\langle \eps \left(\frac{x-1}{x+1}\right)^2\right\rangle_{\hskip-4pt 6}\cup
\left\langle \eps \left(\frac{x-i}{x+i}\right)^2\right\rangle_{\hskip-4pt 6}\\
=&\left\{ 1-x^2,x^2,\frac{1}{1-x^2},\frac{1}{x^2},\frac{x^2-1}{x^2},\frac{x^2}{x^2-1},
\frac{(x+i)^2}{(x-i)^2},\frac{(x+i)^2}{4ix},\frac{(x-i)^2}{(x+i)^2},\frac{(x-i)^2}{-4ix},\frac{-4ix}{(x-i)^2},\frac{4ix}{(x+i)^2},\right.\\
&-x^2,x^2+1,\frac{-1}{x^2},\frac{1}{x^2+1},\frac{x^2+1}{x^2},\frac{x^2}{x^2+1},
\frac{(x+i)^2}{(x-i)^2},\frac{(x+i)^2}{4ix},\frac{(x-i)^2}{(x+i)^2},\frac{(x-i)^2}{-4ix},\frac{-4ix}{(x-i)^2},\frac{4ix}{(x+i)^2},\\
&\left.\frac{1+x^2}{2},\frac{1-x^2}{2},\frac{2}{1-x^2},\frac{2}{x^2+1},\frac{x^2+1}{x^2-1},\frac{x^2-1}{x^2+1},
\frac{2x^2}{x^2+1},\frac{2x^2}{x^2-1},\frac{x^2+1}{2x^2},\frac{x^2-1}{2x^2},\frac{1-x^2}{x^2+1},\frac{x^2+1}{1-x^2}\right\},\\
\ \hskip-2cm
O(x^4):=&
\langle x^4\rangle_6\cup
\left\langle \left(\frac{x+1}{x-1}\right)^4\right\rangle_{\hskip-4pt 6}\cup
\left\langle \left(\frac{x-i}{x+i}\right)^4\right\rangle_{\hskip-4pt 6}\\
=&\left\{ 1-x^4,x^4,\frac{1}{1-x^4},\frac{1}{x^4},\frac{x^4-1}{x^4},\frac{x^4}{x^4-1},
\frac{(x+1)^4}{8x(x^2+1)},\frac{8x(x^2+1)}{(x+1)^4},\frac{-(x-1)^4}{8x(x^2+1)},\frac{-8x(x^2+1)}{(x-1)^4},\right.\\
&\left.\frac{(x+1)^4}{(x-1)^4},\frac{(x-1)^4}{(x+1)^4},\frac{(x+i)^4}{(x-i)^4},\frac{(x-i)^4}{(x+i)^4},\frac{8ix(1-x^2)}{(x-i)^4},\frac{8ix(x^2-1)}{(x+i)^4},
                        \frac{i(x-i)^4}{8x(x^2-1)},\frac{-i(x+i)^4}{8x(x^2-1)}\right\},\\
\ \hskip-2cm
O\left(\frac{2x}{x^2+1}\right):=& \bigcup_{\eps =\pm 1}
\left\langle \frac{2\eps x}{x^2+1}\right\rangle_{\hskip-4pt 6}\cup
\left\langle \frac{2\eps ix}{x^2-1}\right\rangle_{\hskip-4pt 6}\cup
\left\langle \frac{\eps(x^2-1)}{x^2+1}\right\rangle_{\hskip-4pt 6}\\
=
&\left\{ \frac{2x}{x^2+1},\frac{x^2+1}{2x},\frac{(x-1)^2}{-2x},\frac{(x-1)^2}{x^2+1},\frac{-2x}{(x-1)^2},\frac{x^2+1}{(x-1)^2},\right.\\
&\frac{(x+1)^2}{2x},\frac{2x}{(x+1)^2},\frac{-2x}{x^2+1},\frac{x^2+1}{-2x},\frac{(x+1)^2}{x^2+1},\frac{x^2+1}{(x+1)^2}, \\
&\frac{2ix}{x^2-1},\frac{x^2-1}{2ix},\frac{(x-i)^2}{-2ix},\frac{-2ix}{(x-i)^2},\frac{(x-i)^2}{x^2-1},\frac{x^2-1}{(x-i)^2},
\frac{1+x^2}{2},\frac{2}{x^2+1},\frac{1-x^2}{2},\frac{2}{1-x^2},\frac{x^2+1}{x^2-1},\frac{x^2-1}{x^2+1},\\
&\left.\frac{-2ix}{x^2-1},\frac{x^2-1}{-2ix},\frac{(x+i)^2}{2ix},\frac{2ix}{(x+i)^2},\frac{(x+i)^2}{x^2-1},\frac{x^2-1}{(x+i)^2},
\frac{2x^2}{x^2+1},\frac{x^2+1}{2x^2},\frac{2x^2}{x^2-1},\frac{x^2-1}{2x^2},\frac{1-x^2}{x^2+1},\frac{x^2+1}{1-x^2}\right\},\\
\ \hskip-2cm
\left\langle \frac{4x^2}{(x^2+1)^2}\right\rangle_{\hskip-4pt 6}
=&\left\{\frac{(x^2+1)^2}{4x^2},\frac{(x^2-1)^2}{-4x^2},\frac{-4x^2}{(x^2-1)^2},
\frac{4x^2}{(x^2+1)^2},\frac{(x^2+1)^2}{(x^2-1)^2},\frac{(x^2-1)^2}{(x^2+1)^2}\right\},\\
\ \hskip-2cm
O\left(\frac{2x}{x+1}\right):=&
\bigcup_{\eps =\pm 1}
\left\langle \frac{\eps(x-1)}{x+1}\right\rangle_{\hskip-4pt 6}\cup
\left\langle \frac{\eps(x-i)}{x+i}\right\rangle_{\hskip-4pt 6} \cup
\left\langle \frac{\om x}{x+\eps}\right\rangle_{\hskip-4pt 6}\cup
\left\langle \frac{\bom x}{x+\eps}\right\rangle_{\hskip-4pt 6}\\
=
&\left\{\frac{1-x}{2},\frac{2}{1-x},\frac{2}{x+1},\frac{x+1}{2},\frac{x+1}{x-1},\frac{x-1}{x+1},
\frac{\bom x}{x+i},\frac{x+i}{\bom x},\frac{x-1}{\om x},\frac{\om x}{x-1},\frac{i(x+i)}{x-1},\frac{x-1}{i(x+i)},\right.\\
&\frac{-2i}{x-i},\frac{x-i}{-2i},\frac{2i}{x+i},\frac{x+i}{2i},\frac{x+i}{x-i},\frac{x-i}{x+i},
\frac{\bom x}{x-1},\frac{x-1}{\bom x},\frac{x-i}{\om x},\frac{\om x}{x-i},\frac{1-x}{i(x-i)},\frac{i(x-i)}{1-x},\\
&\frac{1-x}{x+1},\frac{x+1}{1-x},\frac{x+1}{2x},\frac{2x}{x+1},\frac{2x}{x-1},\frac{x-1}{2x},
\frac{\bom x}{x+1},\frac{x+1}{\bom x},\frac{x+i}{\om x},\frac{\om x}{x+i},\frac{i(x+1)}{x+i},\frac{-i(x+i)}{x+1},\\
&\left.\frac{2x}{x+i},\frac{x+i}{2x},\frac{x-i}{2x},\frac{2x}{x-i},\frac{i-x}{x+i},\frac{x+i}{i-x},
\frac{\bom x}{x-i},\frac{x-i}{\bom x},\frac{\om x}{x+1},\frac{x+1}{\om x},\frac{i(x-i)}{x+1},\frac{x+1}{i(x-i)}\right\},\\
\ \hskip-3cm
O\left(\frac{x(x-1)}{x^2+1}\right):= &
\left\langle \frac{x(x-1)}{x^2+1}\right\rangle_{\hskip-4pt 6}\cup
\left\langle \frac{x(x+1)}{x^2+1}\right\rangle_{\hskip-4pt 6}\cup
\left\langle \frac{x(x-i)}{x^2-1}\right\rangle_{\hskip-4pt 6}\cup
\left\langle \frac{x(x+i)}{x^2-1}\right\rangle_{\hskip-4pt 6}\\
=
&\left\{\frac{1+x^2}{1-x},\frac{1-x}{1+x^2},\frac{x(x+1)}{x-1},\frac{x-1}{x(x+1)},\frac{x(x+1)}{x^2+1},\frac{x^2+1}{x(x+1)},\right.\\
&\frac{x(1-x)}{x+1},\frac{x+1}{x(1-x)},\frac{x+1}{x^2+1},\frac{x^2+1}{x+1},\frac{x(x-1)}{x^2+1},\frac{x^2+1}{x(x-1)},\\
&\frac{-ix(x+i)}{x-i},\frac{1-x^2}{i(x-i)},\frac{i(x-i)}{1-x^2},\frac{x(x+i)}{x^2-1},\frac{x^2-1}{x(x+i)},\\
&\left. \frac{i(x+i)}{x(i-x)},\frac{x(i-x)}{i(x+i)},\frac{i(x+i)}{x^2-1},\frac{x^2-1}{i(x+i)},\frac{x(x-i)}{x^2-1},\frac{x^2-1}{x(x-i)}\right\},\\
\ \hskip-2cm
O\left(\frac{2(1+i)x}{(x+1)(x+i)}\right):= &
\left\langle \frac{2(1+i)x}{(x+1)(x+i)}\right\rangle_{\hskip-4pt 6}\cup
\left\langle \frac{2(1-i)x}{(x+1)(x-i)}\right\rangle_{\hskip-4pt 6}\\
=
&\left\{\frac{2(1+i)x}{(x+1)(x+i)},\frac{(x+1)(x+i)}{2(1+i)x}, \frac{(x-1)(x-i)}{(x+1)(x+i)},
\frac{(x+1)(x+i)}{(x-1)(x-i)}, \frac{-2(1+i)x}{(x-1)(x-i)},\frac{(x-1)(x-i)}{-2(1+i)x},\right.\\
&\left.\frac{2(1-i)x}{(x+1)(x-i)}, \frac{(x+1)(x-i)}{2(1-i)x}, \frac{(x-1)(x+i)}{(x+1)(x-i)},
\frac{(x+1)(x-i)}{(x-1)(x+i)}, \frac{-2(1-i)x}{(x-1)(x+i)}, \frac{(x-1)(x+i)}{-2(1-i)x}\right\}.\\
\end{align*}

\section*{Appendix B. All $3$-unital functions}
Set $u=\exp(\frac{2\pi i}3)$. In this paper we have proved that there are 84 $3$-unital functions which are given by
\begin{align*}
\vU_3=&\langle x^3\rangle_6\cup  \left\langle\frac{(1-\baru)x}{x-u}\right\rangle_{\hskip-4pt 6} \cup
\bigcup_{\eps\in\gG_3} \langle \eps x\rangle_6  \cup \left\langle \frac{1-u}{\eps x-u}\right\rangle_{\hskip-4pt 6}  \cup\left\langle \frac{1-u\eps x}{1-\eps x}\right\rangle_{\hskip-4pt 6} \cup  \left\langle \frac{-3\eps x}{(\eps x-1)^2}\right\rangle_{\hskip-4pt 6}\\
=& \left\{
x, 1-x, \frac{x}{x-1},\frac{x-1}{x}, \frac{1}{x},\frac{1}{1-x}, ux, 1-ux, \frac{ux-1}{ux},\frac{1}{1-ux},\frac{1}{ux}, \frac{ux}{ux-1},\right.\\
& \baru x, 1-\baru x, \frac1{\baru x}, \frac{1}{1-\baru x}, \frac{\baru x}{\baru x-1}, \frac{\baru x-1}{\baru x}, x^3, 1-x^3, \frac{1}{x^3}, \frac{1}{1-x^3},\frac{x^3}{x^3-1}, \frac{x^3-1}{x^3}, \\
& \frac{1-u}{x-u}, \frac{x-u}{1-u}, \frac{x-1}{u-1}, \frac{u-1}{x-1}, \frac{x-u}{x-1}, \frac{x-1}{x-u},
\frac{1-\baru }{x-\baru }, \frac{x-\baru }{1-\baru }, \frac{x-1}{\baru -1},\frac{\baru -1}{x-1}, \frac{x-\baru }{x-1}, \frac{x-1}{x-\baru }, \\
& \frac{(1-\baru)x}{x-u},\frac{x-u}{(1-\baru)x},\frac{(1-u)x}{x-\baru}, \frac{x-\baru}{(1-u)x}, \frac{x-\baru}{u(x-u)},\frac{u(x-u)}{x-\baru },\\
& \frac{x-\baru }{u-\baru }, \frac{u-\baru }{x-\baru },\frac{\baru -u}{x-u}, \frac{x-u}{\baru -u}, \frac{x-u}{x-\baru }, \frac{x-\baru }{x-u}, \\
& \frac{1-ux}{x-u}, \frac{x-u}{1-ux}, \frac{x-\baru }{u(x-u)}, \frac{u(x-u)}{x-\baru }, \frac{x-1}{\baru -ux},\frac{x-1}{u-\baru x}, \\
& \frac{1-ux}{1-x}, \frac{1-x}{1-ux}, \frac{(1-u)x}{x-1},\frac{x-1}{(1-u)x}, \frac{(1-\baru)x}{x-\baru },\frac{x-\baru }{(1-\baru)x},\\
& \frac{(1-\baru)x}{x-1}, \frac{x-1}{(1-\baru)x}, \frac{1-\baru x}{1-x}, \frac{1-x}{1-\baru x}, \frac{(1-u)x}{x-u}, \frac{x-u}{(1-u)x}, \\
& \frac{-3x}{(x-1)^2}, \frac{(x-1)^2}{-3x}, \frac{(x-1)^2}{x^2+x+1}, \frac{x^2+x+1}{(x-1)^2}, \frac{x^2+x+1}{3x}, \frac{3x}{x^2+x+1}, \\
&\frac{-3\baru x}{(x-\baru)^2},\frac{(x-\baru)^2}{-3\baru x}, \frac{(x-1)(x-u)}{(x-\baru)^2}, \frac{(x-\baru)^2}{(x-1)(x-u)},\frac{3\baru x}{(x-1)(x-u)}, \frac{(x-1)(x-u)}{3\baru x}, \\
&\left.  \frac{-3ux}{(x-u)^2}, \frac{(x-u)^2}{-3ux}, \frac{(x-1)(x-\baru)}{(x-u)^2}, \frac{(x-u)^2}{(x-1)(x-\baru)}, \frac{3ux}{(x-1)(x-\baru)}, \frac{(x-1)(x-\baru)}{3u x} \right\}.
\end{align*}

\end{document}